\documentclass{amsart}
\usepackage{graphicx}
\usepackage{amssymb}
\usepackage{amsmath}
\usepackage{amsthm}
\usepackage{subfigure}
\usepackage{enumerate}
\usepackage{bbm}
\usepackage[]{algorithm2e}

\usepackage[utf8]{inputenc}
\usepackage{amsfonts}
\usepackage{xcolor}
\usepackage{graphicx}

\newtheorem{rem}{Remark}
\newtheorem{lemma}{Lemma}
\newtheorem{thm}{Theorem}
\newtheorem{cor}{Corollary}
\newtheorem{prop}{Proposition}
\newtheorem{defn}{Definition}

\begin{document}

\title[On singularities of the Gauss map components] {On singularities of the Gauss map components of surfaces in ${\mathbb R}^4$}
\date{\today}
\author{W. Domitrz}
\address{Faculty of Mathematics and Information Science, Warsaw University of Technology,  ul. Koszykowa 75, 00-662 Warszawa, Poland.}
\email{wojciech.domitrz@pw.edu.pl}
\author{L. I. Hern\'andez-Mart\'inez}  
\address{Universidad Aut\'onoma de la Ciudad de M\'exico--San Lorenzo Tezonco, (UACM) Prolongaci\'on San Isidro 151, San Lorenzo Tezonco, 09790 Iztapalapa, CDMX, M\'exico.}
\email{lucia.hernandez@uacm.edu.mx}
\author{F. S\'anchez-Bringas}
\address{Universidad Nacional Aut\'onoma de M\'exico (UNAM) Facultad de Ciencias, Ciudad Universitaria, 14210 Coyoacan, CDMX, M\'exico.}
\email{ sanchez@servidor.unam.mx}


\maketitle

\begin{abstract}
The Gauss map of a generic immersion of a smooth, oriented surface into $\mathbb R^4$ is an immersion. But this map takes values on the Grassmanian of oriented 2-planes in $\mathbb R^4$. Since this manifold has a structure of a product of two spheres, the Gauss map has two components that take values on the sphere. We study the singularities of the components of the Gauss map and relate them to the geometric properties of the generic immersion. Moreover, we prove that the singularities are generically stable, and we connect them to  the contact type of the surface and $\mathcal J$-holomorphic curves with respect to an orthogonal complex structure $\mathcal J$ on $\mathbb R^4$.   
Finally, we get some formulas of Gauss-Bonnet type involving the geometry of the singularities of the components with the geometry and  topology of the surface.  
\end{abstract}

\section{Introduction}
The Gauss map singularities of an immersion of a smooth, oriented surface $M$ into $\mathbb R^3$ have been a subject of research because 
they represent a strong relationship between the singularity theory and the geometry of the immersions of a surface. 
Since this map can be locally described as a map from the plane to the plane, the Whitney theory (\cite{Wh}) on singularities of 
differentiable maps is applied.
These singularities arise from the contact of the surface with its tangent plane. The analysis of this contact goes back to classic 
differential geometry of surfaces and has been developed in several works; for instance, \cite{BGM}, \cite{BW}, and \cite{IRST} for an
extensive list of references.
The parabolic set, where the Gaussian curvature vanishes, is the singular set of the map. 
Most interesting properties from the smooth point of view occur near this set. If $M$ is closed,  
the parabolic set of the image of a generic immersion consists of a finite number of 
embedded circles constituted by stable singularities. That is, fold points and isolated cusp points \cite{BW}.
This structure of the parabolic set allows us to find geometric properties of generic immersions of $M$ into $\mathbb R^3$.

If $M$ is immersed into $\mathbb R^4$ the Gauss map $g$ takes values into the Grassmannian of oriented 2-planes in $\mathbb R^4$. This manifold has a structure
of the product $S_1\times S_2$ of 2-dimensional spheres $S_1=S_2=S^2\left(\frac{1}{\sqrt{2}}\right)$ of radius $\frac{1}{\sqrt{2}}$ in $\mathbb R^3$. Then, $g$ has two component-maps from the surface into the sphere. Namely,
$g=(g_1,g_2)=(\pi_1\circ G,\pi_2 \circ G)$, where $\pi_i:S_1 \times S_2 \rightarrow S_i$ is the projection for $i=1,2$.  
Despite several authors studied the Gauss map in this setting, see for instance  \cite{B}, \cite{CS}, \cite{HO1}, and
\cite{We}, an analysis of the singularities of $g_i,\ i=1,2$ has not been developed yet.  
This article aims to settle down the first steps on the study of the singularities of the components of the Gauss map 
and relate them to  
geometric properties of the immersed surface. We begin by characterizing the kernel of the derivative of $g_i,\ i=1,2,$ at each point, as the kernel of a Pfaffian system defined by the connection forms independently of the coordinate chart (Theorem \ref{kernel}).   
In \cite{HO1} (see also \cite{We}), the authors proved that
the singular set of the component $g_i=\pi_i\circ G$ for $i=1$ (or $i=2$) is the zero locus of the function defined as the sum (or the difference, respectively) of the Gaussian and normal curvatures. In the present article, following the ideas from \cite{B}, we provide an alternate proof of this result
(Corollary \ref{c1}) which is rather simple and more convenient to the approach of this research.
 
The Gauss map of a smooth, oriented surface immersed into $\mathbb R^4$  might not be regular, i.e., it might not be an immersion at some points.  
The lack of regularity imposes strong restrictions on the surface. 
Specifically, these non-regular points are characterized as the points where 
the Gaussian curvature and the discriminant of the second fundamental form vanish simultaneously, see \cite{L} for other equivalent conditions.
We prove that for a generic immersion of $M$ into $\mathbb R^4$, the Gauss map is regular (Proposition \ref{codim3}).
Further on, we define three  
properties (G1), (G2) and (G3) for immersions of $M$ into $\mathbb R^4$ that determine the 
structure of the singular set of the components of the Gauss map. Namely, for immersions having these properties
the singular set of $g_i,\ i=1,2$ consists of a finite number of smooth embedded circles constituted by stable singularities. That is, fold points and isolated cusp points. 
Moreover, the corresponding singular set's image under $g_i\ i=1,2$ has stable intersections.
Then, we prove the genericity theorem. Namely, immersions having these properties are generic (Theorem \ref{generic}).
Contrasting with the case of the Gauss map of surfaces in $\mathbb R^3$,
the singularities of the Gauss map components of surfaces in $\mathbb R^4$ arise from the contact of the surface with a $\mathcal J$-holomorphic curve with respect to an orthogonal complex structure $\mathcal J$  in 
$\mathbb R^4$ that better approaches the surface locally. We introduce the study of this type of contact in the present article. So, we provide a classification of regular and generic singular points of these components employing the normal forms of jets of contact maps  (Theorem \ref{NormalForms}). 
Further on, we apply certain theorems proved in \cite{Q}, \cite{SUY} and \cite{SUY3} to get formulas of Gauss-Bonnet type relating the geometry and topology of the surface to singularities of the components of $g$ 
(Theorem \ref{G-B} and Corollary \ref{G-B-1}). 
The article is organized as follows. In section 2, we present some preliminary results including the description of the second order invariants of an immersion of a surface into $\mathbb R^4$ in terms of the connection forms.
Moreover, we describe the model 
of the oriented Grassmannian of 2-planes in $\mathbb R^4$ that will be employed in the research. Section 3 is devoted to studying the 
singular sets of the components of the Gauss map; we prove Theorem \ref{kernel} and Corollary \ref{c1}. In section 4, we define the genericity conditions $(G_1)$, $(G_2)$ and $(G_3)$ and we discuss their meaning concerning the structure of the singular sets of $g_i,\ i=1,2$.   
In section 5, we show the local description of the components of the Gauss map. 
We use this fundamental tool in the proof of the genericity theorem.  
In section 6, we deal with the generic properties of the Gauss map and its components. Specifically, we prove Proposition \ref{codim3} and  Theorem \ref{generic}. Section 7 aims to prove the classification theorem of regular and generic singular point.  
In section 8, we demonstrate Theorem \ref{G-B} and Corollary \ref{G-B-1} that provide the Gauss-Bonnet type formulas.
We conclude by analyzing the example of the double torus mentioned above. 
Finally, we point out that the research presented in this article has a relevant extension by using the approach presented in
\cite{DZ} for the class of oriented surfaces with boundary.  
      
 \section{Preliminaries}

\subsection{The second-order invariants}
Let $M$ be a surface immersed into $\mathbb R^4$. We denote by ${\mathcal X}(M)$ and ${\mathcal X}(M)^{\perp}$ the
spaces of tangent and normal smooth vector fields on $M$, respectively. 
The second fundamental form of $M$ is the tensor field  
defined as
\[
\begin{array}{ccccl}
II  & : & {\mathcal X}(M) \times {\mathcal X}(M)  \longrightarrow 
{\mathcal X}(M)^{\perp}, \ \ \ 
II(X, Y)= (dY(X))^{\perp},
\end{array}
\]
where
$Z^{\perp}$ denotes the orthogonal projection   on ${\mathcal X}(M)^{\perp}$ of a smooth vector field $Z:M\rightarrow \mathbb R^4$. 
This tensor field determines at each point $p \in M$ a quadratic map:
\[
\begin{array}{ccccl}
II_p   :  T_pM   \longrightarrow & T_pM^{\perp},\ \ \  II_p (X) = II(X,X)_p
\end{array}
\]
whose invariants under the isometries of $T_pM$ and  $T_pM^{\perp}$, respectively,  
are 
the second-order invariants of $M$. They are known as the norm of the mean curvature vector,
the Gaussian curvature, the normal curvature and the discriminant of the second fundamental form.  We denote them by 
\begin{eqnarray}\label{invariants}
|H|^2,\ K,\ K^N,\ \Delta,
\end{eqnarray}
 respectively; see \cite{L} and 
\cite{BS}.

Let $U \subset \mathbb R^2$ be an open set parameterized by $(u,v)$. We suppose that ${\bf x}:U \rightarrow M$ is a local parameterization of $M$ endowed with 
a Darboux frame  $(e_1, e_2, e_3, e_4)$. Namely, a 
positive orthonormal frame of $\mathbb R^4$ on  ${\bf x}(U)$ such that, for every $(u,v)\in U$, $(e_1(u,v), e_2(u,v))$ 
is a positive orthonormal basis of $T_{q}M$, and $(e_3(u,v), e_4(u,v))$ is an orthonormal basis of $T_qM^{\perp}$, where $\ q={\bf x}(u,v)$. 
We denote by $\theta_i=\left\langle d{\bf x},e_i\right\rangle $ for $i=1,2,3,4$ the dual forms of the frame. 
Notice that $\theta_3=\theta_4=0$. Moreover,
the connection $1$-forms $\omega_i^j$ for $i,j=1,\cdots,4$ are defined in the following way 
\begin{equation}\label{connection}
de_i(X)=\sum_{j=1}^4\omega_i^j (X) e_j,
\end{equation}
where $X$ is a tangent vector field on $M$. 
Thus, 
the curvatures $K$ and $K^N$ satisfy 

\begin{equation}\label{streq}
  \begin{array}{ccl}
d\omega_1^2 = -K\ \theta_1 \wedge \theta_2,\ \ \ 
d\omega_3^4 = -K^N \theta_1 \wedge \theta_2,
  \end{array}
\end{equation}
where $\theta_1 \wedge \theta_2$ is the area form defined by the tangent frame.

The second fundamental form is expressed as

\begin{equation}\label{streq}
  \begin{array}{ccl}
II(X)&=&(\omega_1^3(e_1)(\theta_1)^2+(\omega_2^3(e_1)+\omega_1^3(e_2))\theta_1 \theta_2+\omega_2^3 (e_2)(\theta_2)^2)e_3 \\
&+&(\omega_1^4(e_1)(\theta_1)^2+(\omega_2^4(e_1)+\omega_1^4(e_2))\theta_1 \theta_2+\omega_2^4(e_2)(\theta_2)^2)e_4.
  \end{array}
\end{equation}

\subsection{The oriented Grassmannian $Gr^+(2,4)$}

We denote by $Gr^+(2,4)$ the space of all oriented $2$-dimensional subspaces (planes) of $\mathbb R^4$. Let $(e_1,e_2,e_3,e_4)$ be a positive orthonormal  basis for $\mathbb R^4$, i.e., $\left\langle e_i,e_j\right\rangle =\delta_{ij}$ for $i,j=1,2,3,4$ and $e_1\wedge e_2\wedge e_3 \wedge e_4$ is the orientation of $\mathbb R^4$. We use 
Pl$\ddot{\rm u}$cker coordinates in $Gr^+(2,4)$ (see \cite{CS}). 
Let $\Lambda^2 \mathbb R^4$ be the space of bivectors of $\mathbb R^4$. Thus,
if $V$ is an oriented plane of $\mathbb R^4$ and $(v_1,v_2)$ is a positive orthonormal basis of $V$, the bivector 
\begin{equation}\label{Plucker}
v_1\wedge v_2=\alpha_{12}e_1\wedge e_2+\alpha_{23}e_2\wedge e_3+\alpha_{31}e_3\wedge e_1+\alpha_{34}e_3\wedge e_4+\alpha_{14}e_1\wedge e_4+\alpha_{24}e_2\wedge e_4,
\end{equation}
in the basis  $e_i\wedge e_j$ used in this expression,
represents $V$ and
the coordinates $\alpha_{ij}$ are called Pl$\ddot{\rm u}$cker coordinates. They are independent of the choice of a positive orthonormal basis of $V$ and satisfy the following equations

\begin{equation}\label{sphere}
  \begin{array}{ccl}
\alpha_{12}\alpha_{34}+\alpha_{23}\alpha_{14}+\alpha_{31}\alpha_{24}&=&0,\\
\alpha_{12}^2+\alpha_{34}^2+\alpha_{23}^2+\alpha_{14}^2+\alpha_{31}^2+\alpha_{24}^2&=&1. 
  \end{array}
\end{equation}

The first equation of system (\ref{sphere}) follows from $(v_1\wedge v_2)\wedge (v_1\wedge v_2)=0$. Moreover, we extend the dot product on $\mathbb R^4$ to $\Lambda^2 \mathbb R^4$ by the standard formula 
\begin{equation}
\left\langle u_1\wedge u_2,w_1\wedge w_2\right\rangle =\det\left(\begin{array}{cc}\left\langle u_1,w_1\right\rangle  & \left\langle u_1,w_2\right\rangle \\
\left\langle u_2,w_1\right\rangle &\left\langle u_2,w_2\right\rangle \\
\end{array}\right),
\end{equation}
 for $u_1,u_2,w_1,w_2 \in \mathbb R^4$. Then, the second equation of (\ref{sphere}) follows from 
 $$
 \left\langle v_1\wedge v_2,v_1\wedge v_2\right\rangle =1.
 $$
  Straightforward computations imply that if the coordinates $\alpha_{ij}$ of $w\in \Lambda^2\mathbb R^4$ satisfy system (\ref{sphere}), then $w=w_1\wedge w_2$ for $w_1, w_2 \in \mathbb R^4$, such that  
$\left\langle w_i,w_j\right\rangle =\delta_{ij}$ for $i,j=1,2$. 

The following bivectors
\begin{eqnarray}\label{x1}
\left\lbrace 	
\begin{array}{cc}
	x_1=\frac{1}{\sqrt{2}}(e_1\wedge e_2 + e_3\wedge e_4), & y_1=\frac{1}{\sqrt{2}}(e_1\wedge e_2 - e_3\wedge e_4), \label{x1} \\ 
x_2=\frac{1}{\sqrt{2}}(e_2\wedge e_3 + e_1\wedge e_4), & y_2=\frac{1}{\sqrt{2}}(e_2\wedge e_3 - e_1\wedge e_4), \label{x2}\\
x_3=\frac{1}{\sqrt{2}}(e_3\wedge e_1 + e_2\wedge e_4), & y_3=\frac{1}{\sqrt{2}}(e_3\wedge e_1 - e_2\wedge e_4)  \label{x3} 
\end{array} \right.
\end{eqnarray}
form an orthonormal basis of $\Lambda^2 \mathbb R^4$. System (\ref{sphere}) implies that  coordinates 
of $v_1\wedge v_2=\sum_{i=1}^3 \beta_i x_i+\gamma_i y_i$ satisfy the following equations
\begin{equation}\label{2-sphere}
\sum_{i=1}^3\beta_i^2=\frac{1}{2}, \  \  \ \ \  \sum_{i=1}^3\gamma_i^2=\frac{1}{2}.
\end{equation}
So, $Gr^+(2,4)$ is the product of two $2$-dimensional spheres $S_{1}\times S_{2}$ of radius $\frac{1}{\sqrt{2}}$ in $\mathbb R^3$.

\section{The singular sets of Gauss map components.}

Let $(e_1, e_2, e_3, e_4)$  be a local Darboux frame on $M$. Then Cartan's first  structural equations, $d\theta=\omega\wedge \theta$  imply for $j=3,4$ 
\begin{equation*}
0=d\theta_j=\omega_j^1\wedge \theta_1+\omega_j^2 \wedge \theta_2, 
\end{equation*}
which means that we have
\begin{equation}\label{w1j2=w2j1}
\omega_1^j(e_2)=\omega_2^j(e_1) \ \text{for} \ j=3,4. 
\end{equation}
By Cartan's second structural equations $d\omega=\omega\wedge \omega$, we have
\begin{eqnarray*}
d\omega_1^2&=\omega_1^3\wedge\omega_3^2+\omega_1^4\wedge\omega_4^2&=-\omega_1^3\wedge\omega_2^3-\omega_1^4\wedge\omega_2^4,\\
d\omega_3^4&=\omega_3^1\wedge\omega_1^4+\omega_3^2\wedge\omega_2^4&=-\omega_1^3\wedge\omega_1^4-\omega_2^3\wedge\omega_2^4.
\end{eqnarray*}
Therefore, 

\begin{equation}\label{domega12}
  d\omega_1^2+d\omega_3^4=-(\omega_1^3-\omega_2^4)\wedge(\omega_1^4+\omega_2^3), \
d\omega_3^4-d\omega_1^2=(\omega_1^4-\omega_2^3)\wedge(\omega_2^4+\omega_1^3). 
 \end{equation}
So, by evaluating system (\ref{streq}) on the basis $(e_1,e_2)$, we obtain that 
\begin{eqnarray}\label{detC}
K+K^N&=&\det\left(\begin{array}{cc}
\omega_1^3(e_1)-\omega_2^4(e_1)&\omega_1^3(e_2)-\omega_2^4(e_2)\\
\omega_1^4(e_1)+\omega_2^3(e_1)&\omega_1^4(e_2)+\omega_2^3(e_2)\\
\end{array}\right),\\
K-K^N&=&\det\left(\begin{array}{cc}
\omega_1^4(e_1)-\omega_2^3(e_1)&\omega_1^4(e_2)-\omega_2^3(e_2)\\
\omega_2^4(e_1)+\omega_1^3(e_1)&\omega_2^4(e_2)+\omega_1^3(e_2)\\

\end{array}\right).
\end{eqnarray}
These equations imply the following equivalences.
\begin{prop}
Let $p$ be a point in $M$. 

$(K+K^N)(p)=0$ if and only if $1$-forms $\omega_1^3-\omega_2^4$, $\omega_1^4+\omega_2^3$ are linearly dependent at $p$.

$(K-K^N)(p)=0$ if and only if $1$-forms $\omega_1^3+\omega_2^4$, $\omega_1^4-\omega_2^3$ are linearly dependent at $p$.
\end{prop}

\begin{lemma}
The kernels of Pfaffian systems $\{\omega_1^3-\omega_2^4=0,\, \omega_1^4+\omega_2^3=0\}$ and $\{\omega_1^4-\omega_2^3=0,\, \omega_2^4+\omega_1^3=0\}$ on an open subset $U$ of $M$ do not depend on the choice of a local Darboux frame $(e_1, e_2, e_3, e_4)$ on $U$. 
\end{lemma}
\begin{proof}
Any local Darboux frame on $U$ has the following form $(Ae_1, Ae_2, Ae_3, Ae_4)$, where
\begin{equation}\label{transformation}
A=\left(\begin{array}{cccc}
\cos\alpha&-\sin\alpha&0&0\\
\sin\alpha&\cos\alpha&0&0\\
0&0&\cos\beta&-\sin\beta\\
0&0&\sin\beta&\cos\beta\\
\end{array}\right).
\end{equation}
If $\omega$ is a connection form for $(e_1,e_2,e_3,e_4)$ then $\tilde{\omega}=dA A^{-1}+A\omega A^{-1}$ is a connection form for $(Ae_1, Ae_2, Ae_3, Ae_4)$.

Therefore, we get
\begin{eqnarray*}
\tilde{\omega}_1^3-\tilde{\omega}_2^4&=&\cos(\alpha + \beta) (\omega_1^3 - \omega_2^4) - \sin(\alpha + \beta)(\omega_1^4 + \omega_2^3), \\
\tilde{\omega}_1^4+\tilde{\omega}_2^3&=&\sin(\alpha + \beta) (\omega_1^3 - \omega_2^4) + \cos(\alpha + \beta)(\omega_1^4 + \omega_2^3), \\
\tilde{\omega}_1^3+\tilde{\omega}_2^4&=&\cos(\beta-\alpha)(\omega_1^3 + \omega_2^4)-\sin(\beta-\alpha) (\omega_1^4 - \omega_2^3),\\
\tilde{\omega}_1^4-\tilde{\omega}_2^3&=& \sin(\beta-\alpha)(\omega_1^3 + \omega_2^4)+\cos(\beta-\alpha) (\omega_1^4 - \omega_2^3), 
\end{eqnarray*}
which finish the proof.
\end{proof}

Let $(e_1,e_2,e_3,e_4)$ be a local Darboux frame on $U$. The Gauss map (in the Pl$\ddot{\rm u}$cker coordinates) has the following form
\begin{eqnarray}\label{gauss_in_Plucker}
g: M \rightarrow  Gr^+(2,4),\ \ 
   p \mapsto e_1(p)\wedge e_2(p)  
\end{eqnarray} 
So, we express the differential of the Gauss map as 
\begin{equation}\label{dgauss}
dg=d(e_1\wedge e_2)=de_1\wedge e_2-de_2\wedge e_1=\omega_1^3 e_3\wedge e_2+\omega_1^4 e_4\wedge e_2+\omega_2^3 e_1\wedge e_3+\omega_2^4 e_1\wedge
e_4.
\end{equation} 
Therefore, we conclude:
\begin{prop}
The Gauss map is an immersion at $q\in M$ if and only if the dimension of the space $\text{span}\{ \omega_1^3|_q, \omega_1^4|_q, \omega_2^3|_q,\omega_2^4|_q\}$ is greater than $1$. 
\end{prop}
\begin{proof} By (\ref{dgauss}), a tangent vector $v\in T_qM$ is in the kernel of the Gauss map if and only if $\omega_1^3|_q(v)=\omega_1^4|_q(v)=\omega_2^3|_q(v)=\omega_2^4|_q(v)=0$.  If $v\ne 0$ the dimension of the space $\text{span}\{ \omega_1^3|_q, \omega_1^4|_q, \omega_2^3|_q,\omega_2^4|_q\}$ is not greater than $1$.
\end{proof}
 
Lemma $1$ implies that
the dimension of the space $\text{span}\{ \omega_1^3|_q, \omega_1^4|_q, \omega_2^3|_q,\omega_2^4|_q\}$ does not depend on the choice of a local Darboux frame. Thus, the dimension of $\text{span}\{ \omega_1^3|_q, \omega_1^4|_q, \omega_2^3|_q,\omega_2^4|_q\}$ is $0$ if and only if $\omega_i^j|_q=0$ for $i=1,2$ and $j=3,4$. Let us assume that $\omega_1^3|_q\ne 0$. Then $\dim\text{span}\{ \omega_1^3|_q, \omega_1^4|_q, \omega_2^3|_q,\omega_2^4|_q\}=1$ if and only if 
\begin{equation}\label{dim1}
\omega_1^3|_q\wedge \omega_1^4|_q=\omega_1^3|_q\wedge \omega_2^3|_q=\omega_1^3|_q\wedge \omega_3^4|_q=0.
\end{equation}
We use (\ref{dim1}) to prove that the Gauss map is an immersion for a generic surface immersed into $\mathbb R^4$  (see Proposition \ref{codim3}).

We describe the Gauss map (\ref{gauss_in_Plucker}) as 
\begin{equation}\label{g_in_Plucker}
g(p)=e_1(p)\wedge e_2(p)=\frac{1}{\sqrt{2}}x_1(p)+\frac{1}{\sqrt{2}}y_1(p),
\end{equation}
where $(x_1,x_2,x_3,y_1,y_2,y_3)$ is an orthonormal basis on $\Lambda^2 \mathbb R^4$ given by (\ref{x1}).
Equation (\ref{2-sphere}) implies that 
$$\frac{1}{\sqrt{2}}x_1\in S_{1}\subset \text{span}(x_1,x_2,x_3),  \ \frac{1}{\sqrt{2}}y_1\in S_{2}\subset \text{span}(y_1,y_2,y_3).$$
Therefore, the components of the Gauss map $g=(g_1,g_2):M\rightarrow S_{1}\times S_{2}$ have the following form
\begin{equation}\label{g_1}
  \begin{array}{ccl}
g_1(p)=\frac{1}{\sqrt{2}}x_1(p)=\frac{1}{2}(e_1(p)\wedge e_2(p)+e_3(p)\wedge e_4(p)) , \\
 g_2(p)= \frac{1}{\sqrt{2}}y_1(p)=\frac{1}{2}(e_1(p)\wedge e_2(p)-e_3(p)\wedge e_4(p)). 
  \end{array}
\end{equation}
It is straightforward to prove that $g_1$ and $g_2$ do not depend on the choice of a local Darboux frame. 
By deriving them, we get

\begin{equation}\label{dg_1}
  \begin{array}{ccl}
dg_1&=&\frac{1}{\sqrt{2}}(\omega_2^4-\omega_1^3)x_2+\frac{1}{\sqrt{2}}(-\omega_1^4-\omega_2^3)x_3,\\
dg_2&=&\frac{1}{\sqrt{2}}(-\omega_2^4-\omega_1^3)y_2+\frac{1}{\sqrt{2}}(\omega_1^4-\omega_2^3)y_3.
  \end{array}
\end{equation}

Therefore, we state the following.
\begin{thm} \label{kernel}
The kernel of  $dg_1$ coincides with the kernel of the Pfaffian system $\{\omega_1^3-\omega_2^4=0,\, \omega_1^4+\omega_2^3=0\}$. Moreover, 
the kernel of $dg_2$ coincides with the kernel of the Pfaffian system $\{\omega_2^4+\omega_1^3=0,\, \omega_1^4-\omega_2^3=0\}$.
\end{thm}

Let $\sigma$ be the area form on $S_{i}$, $i=1,2$ and $dA$ be the area form on $M$. 
Given a smooth map $f:M\rightarrow S_{i}$, we denote by $J(f)$ the {\it Jacobian} of f. Namely, the smooth function defined by the following equation $f^{\ast}\sigma=J(f)dA$. Formulas (\ref{dg_1}) imply the following classical result.

\begin{thm}[\cite{B}, Proposition 4.5 in  \cite{HO1}, \cite{We}] \label{Jacobian}
Let $J(g_i)$ be the Jacobian of $g_i$ for $i=1,2$ and $K$, $K^N$ the Gaussian and normal curvature, respectively. Then
\begin{eqnarray}
J(g_1)=\frac{1}{2}(K+K^N), \ \ \  
J(g_2)=\frac{1}{2}(K-K^N).  \label{jacg1}
\end{eqnarray}
\end{thm}

\begin{proof}  
Equation (\ref{2-sphere}) implies that $x_2(p),\ x_3(p)$ span
the tangent space to $S_{1}$ at the point $g_1(p)=\frac{x_1(p)}{\sqrt{2}}$, and $y_2(p),\ y_3(p)$ span the tangent space to $S_{2}$ at $g_2(p)=\frac{y_1(p)}{\sqrt{2}}$. 
Hence, 
we apply system (\ref{dg_1}) to get 
$$g_1^{\ast}\sigma=\frac{1}{\sqrt{2}}(\omega_2^4-\omega_1^3)\wedge\frac{1}{\sqrt{2}}(-\omega_1^4-\omega_2^3)=\frac{1}{2}(\omega_1^3-\omega_2^4)\wedge(\omega_1^4+\omega_2^3),$$
$$g_2^{\ast}\sigma=\frac{1}{\sqrt{2}}(-\omega_2^4-\omega_1^3)\wedge\frac{1}{\sqrt{2}}(\omega_1^4-\omega_2^3)=\frac{1}{2}(\omega_1^4-\omega_2^3)\wedge(\omega_2^4+\omega_1^3).$$
Thus, system (\ref{domega12}) implies the result.
\end{proof}
Consequently, we state the following.
\begin{cor}\label{c1}
The singular set of the component $g_i,\ i=1,2$ of the Gauss map has the following form 
$$
\Sigma_{g_i}=\{p\in M| (K+(-1)^{i+1}K^N)(p)=0\}.
$$ 
\end{cor}

If the Gauss map $g=(g_1,g_2)$ of a surface $M \subset \mathbb R^4$ is not an immersion at $p$, the kernel of $dg_p$ has dimension greater than zero, so $g_1$ and $g_2$ are singular at $p$ since they are the projections of $g$ onto the spheres.   

The converse of this assertion is not valid. Corollary \ref{c1} implies that
the set of points where $g_1$ and $g_2$ are singular simultaneously is the intersection of the set of flat points, where $K=0$, with 
the set of semiumbilic points, where $K^N=0$.

\section{Generic properties}\label{properties}

Let $\xi(t)$ be a regular local parameterization of a smooth curve $C$  
of $M$ such that $\xi(0)=p$. Let $V$ be a smooth unit vector field along $\xi$, i.e., $V(t)\in T_{\xi(t)}M$ and $<V(t),V(t)>=1$. Then $V$ is {\it $1$-tangent} to $C$ at $p=\xi(0)$ if  $\xi'(0)$ and $V(0)$ are linearly dependent and $\frac{d}{dt}|_{t=0}\left(\xi'(t)\wedge V(t)\right)\ne 0.$ 
A Pfaffian system of codimension 1 along $\xi$ is $1$-tangent to $C$ at $p=\xi(0)$ if the unit generator of its kernel is 
$1$-tangent to $C$ at $p=\xi(0)$.  

Let $C^{\infty}(M,N)$ be the set of all $C^{\infty}$ maps from a compact manifold $M$ to an arbitrary manifold $N$ with the Whitney $C^{\infty}$ topology i.e., the topology of uniform convergence of each $k$-jet ($k=0,1,2,3,\cdots$). Let $I$ be an open subset of $C^{\infty}(M,N)$. We call a property of maps in $I$ { \it generic} if the subset of maps in $I$ with this property is open and dense in $I$. 

The map $f\in C^{\infty}(M,N)$ is {\it stable} if there exists an open neighborhood $U$ of $f$ in $C^{\infty}(M,N)$, such that 
for any $h$ in $U$, $h$ is conjugated to $f$. That is, there exist diffeomorphisms $\Phi:M\rightarrow M$ and $\Psi: N \rightarrow N$, that satisfy $h=\Psi\circ f\circ \Phi$.

In this research $M$ is a smooth, closed, oriented surface and $N=\mathbb R^4$. We mainly study the set $I(M,\mathbb R^4)$ of all immersions of $M$ into $\mathbb R^4$, which is an open subset of $C^{\infty}(M, \mathbb R^4)$.

Let $(e_1,e_2,e_3,e_4)$ be a local Darboux frame on a closed, oriented 
surface $M$ immersed into $\mathbb R^4$ and $\omega_i^j=<de_i,e_j>$.

We define the following properties of immersions in $I(M,\mathbb R^4)$.
\begin{itemize}
\item[$(G_1)$] The $1$-form $d(K\pm K^N)$ never vanishes on the set $\{K\pm K^N=0\}$. 
\item[$(G_2)$] Property $(G_1)$ holds, and the kernel of the Pfaffian system $\{\omega_1^3\mp \omega_2^4=0, \omega_1^4\pm \omega_2^3=0\}$ is transverse to the curve  $\{K\pm K^N=0\}$ except at a finite number of points where it is $1$-tangent.  
\item[$(G_3)$] Properties $(G_1)$ and $(G_2)$ hold and the bivector $\frac{1}{2}(e_1(p)\wedge e_2(p)\pm e_3(p)\wedge e_4(p))$ at any point $p$ such the kernel of the Pfaffian system $\{\omega_1^3\mp \omega_2^3=0, \omega_1^4\pm \omega_2^3=0\}$ is $1$-tangent to the curve  $\{K\pm K^N=0\}$ does not coincide with  $e_1(q)\wedge e_2(q)\pm e_3(q)\wedge e_4(q)$  at any other point $q$ of the curve $\{K\pm K^N=0\}$. There are, at most, finitely many pairs (and no triplets, etc.) of points of the curve  $\{K\pm K^N=0\}$ at which values of $\frac{1}{2}(e_1\wedge e_2\pm e_3\wedge e_4)$ coincide and  the image of $\{K\pm K^N=0\}$ by the map 
$\frac{1}{2}(e_1\wedge e_2\pm e_3\wedge e_4)$ are transverse at these pairs. 
\end{itemize}

Property $(G_1)$ implies that the set $\{K-(-1)^i K^N=0\}$ is a smooth curve embedded in 
$M$ for $i=1,2$. Since $M$ is closed, it is a finite disjoint collection of smoothly embedded circles. By Theorem \ref{Jacobian}, the set $\{K-(-1)^i K^N=0\}$ is the set singular points of $g_i$ for $i=1,2$. 

The generic local rank one singular points are characterized by having normal forms given 
in some pair of local coordinates at $p$ and $f(p)$ by
\begin{eqnarray*}
(u,v)\mapsto (u^2, v)\ \ {\rm or}\ \  (u,v) \mapsto(uv\pm v^3,u).
\end{eqnarray*}
In the first case, $p$ is called a {\it fold point}; in the second, $p$ is called a {\it cusp point}.

By Theorem \ref{kernel}, the kernel of the Pfaffian system $\{\omega_1^3+(-1)^i \omega_2^4=0, \omega_1^4-(-1)^i \omega_2^3=0\}$ coincides with the kernel of $dg_i$ for $i=1,2$. The kernel of $dg_i$ is transverse to the regular curve $\Sigma_{g_i}=\{K-(-1)^i K^N=0\}$ of singular points of $g_i$ at $p$ if and only if $p$ is a fold point of
$g_i$. The kernel of $dg_i$ is $1$-tangent to $\Sigma_{g_i}$ at $p$ if and only if $p$ is a cusp point of $g_i$. 
Thus, property $(G_2)$ implies that the singular points of $g_i$ are fold points except by a finite number of cusp points.

By observing that $g_i=\frac{1}{2}(e_1\wedge e_2-(-1)^i e_3\wedge e_4)$ (see (\ref{g_1})), we conclude that property $(G_3)$ means that the set of singular values of $g_i$ does not intersect itself in the image of the cusp point; there are, at most, finitely many pairs (and no triple, etc.) of images of fold points at which the set of singular values intersects transversally. 

By the very definition (see \cite{Wh}), the map $f$, from $M$  
into another surface, is good if and only if  
the Jacobian and the gradient of the 
Jacobian do not vanish simultaneously.
Moreover, $f$ is excellent if and only if it is good and the singular set consists of fold points or isolated cusp points. 
Thus, an immersion ${\bf x} \in I(M,\mathbb R^4)$ satifies $(G_1)$ if and only if the components of the Gauss map are good. Moreover,
${\bf x}$ satisfies $(G_2)$ if and only if it is excellent.
We will prove that properties $(G_1)$, $(G_2)$ and $(G_3)$ are generic in Section \ref{sing}.

These properties are generalizations to our setting of 
those provided in \cite{BW} for immersions of $M$ into $\mathbb R^3$.

\section{Local description  of the Gauss map components.}

We provide a local description of the Gauss map of a surface $M$ immersed into $\mathbb R^4$ near a point $p$. 
To do so, we translate $p$ to the origin and parameterize $M$  
as the graph of a function in an open neighborhood $U \subset \mathbb R^2$ of the origin, with parameters $(u,v)$, in the following way

\begin{equation*}\label{param}
{\bf x}(u,v)=(u,v,a(u,v),b(u,v)),
\end{equation*}
where $a$ and $b$ are local smooth functions. For simplicity, we introduce the following notation 
\begin{eqnarray*}
	a_{ij}=\frac{\partial^{i+j} a}{\partial u^i \partial v^j}(u,v)\ \  {\rm and}\ \ b_{ij}=\frac{\partial^{i+j} b}{\partial u^i \partial v^j}(u,v), 
\end{eqnarray*}
where $i,j \in \mathbb N \cup \{0 \}$.
If we assume that $a_{ij}$ and $b_{ij}$ vanish at the origin for $i+j=1$, the tangent plane at this point is the $uv$-plane, and we say that $M$ is parameterized locally in {\it Monge form} at the origin.   
The tangent frame on $M$
induced by the parameterization is 
$${\bf x}_u=\left(1,0,a_{10},b_{10}\right) \ \ {\rm and}\ \  
{\bf x}_v=\left(0,1,a_{01},b_{01}\right).$$  
Then, we define the following orthonormal tangent frame with positive orientation

\begin{equation*}\label{tangentframe}
e_1
=\frac{1}{\sqrt{1+ a_{10}^2+b_{10}^2}}\left(\begin{array}{c}1 
\\ 0 
\\ a_{10}  
\\  b_{10}  
\end{array}\right), 
\end{equation*}
\begin{equation*}\label{tangentframe1}
e_2=\frac{{\bf x}_v-\left({\bf x}_v\cdot e_1\right)e_1}{\left|\left({\bf x}_v-{\bf x}_v\cdot e_1\right)e_1\right|}=
\frac{1}{h\sqrt{1+ a_{10}^2+b_{10}^2}}\left(\begin{array}{c}-\left(a_{10}a_{01}+b_{10}b_{01}\right) 
\\ 1+a_{10}^2+b_{10}^2 
\\ a_{01}-a_{10}b_{10}b_{01}+a_{01}b_{10}^2
\\   b_{01}-a_{10}a_{01}b_{10}+b_{01}a_{10}^2  
\end{array}\right),
\end{equation*}
where
\begin{equation}\label{h}
h=\sqrt{1+ a_{10}^2+a_{01}^2+b_{10}^2+b_{01}^2+\left(a_{01} b_{10}-a_{10}b_{01}\right)^2}.
\end{equation}

\noindent Now, 
we complete the tangent frame $\{e_1, e_2\}$   
with the normal frame defined as 
\begin{eqnarray*}\label{E3}
e_3=\frac{1}{\sqrt{1+a_{10}^2 +a_{01}^2}}\left(\begin{array}{c}-a_{10} \\ -a_{01}\\ 1\\ 0\end{array}\right),
\end{eqnarray*}
\begin{eqnarray*}\label{E4}
e_4=\frac{1}{h \sqrt{1+a_{10}^2 +a_{01}^2} }\left(\begin{array}{c}a_{01} a_{10}b_{01}-b_{10}(1+a_{01}^2) \\ a_{01}a_{10}b_{10}-b_{01}(1+a_{10}^2) \\ -a_{01}b_{01}-a_{10}b_{10} \\ 1+a_{10}^2 +a_{01}^2 \end{array}\right),
\end{eqnarray*} 
where $h$ is given by (\ref{h}). We observe that  $e_1\wedge e_2 \wedge e_3 \wedge e_4$ is the positive orientation of $\mathbb R^4$. So, the frame
\begin{eqnarray}\label{darbouxframe}
(e_1,e_2,e_3,e_4)
\end{eqnarray}
is a Darboux frame on $M$.   
Moreover, for the standard positive orthonormal basis  $E_1=(1,0,0,0)$, $E_2=(0,1,0,0)$, $E_3=(0,0,1,0)$, $E_4=(0,0,0,1)$ of $\mathbb R^4$, we define the following orthonormal basis of $\Lambda^2 \mathbb R^4$ 

\begin{eqnarray}\label{X1}
\left\lbrace 	
\begin{array}{cc}
	X_1=\frac{1}{\sqrt{2}}(E_1\wedge E_2 + E_3\wedge E_4), & Y_1=\frac{1}{\sqrt{2}}(E_1\wedge E_2 - E_3\wedge E_4), \\ 
	X_2=\frac{1}{\sqrt{2}}(E_2\wedge E_3 + E_1\wedge E_4), & Y_2=\frac{1}{\sqrt{2}}(E_2\wedge E_3 - E_1\wedge E_4),\\ 
	X_3=\frac{1}{\sqrt{2}}(E_3\wedge E_1 + E_2\wedge E_4), & Y_3=\frac{1}{\sqrt{2}}(E_3\wedge E_1 - E_2\wedge E_4).
	\end{array} \right.
\end{eqnarray}
Then, 
$g_1=\frac{x_1}{\sqrt{2}}=\frac{1}{\sqrt{2}}(e_1\wedge e_2+e_3\wedge e_4)$ and $g_2=\frac{y_1}{\sqrt{2}}=\frac{1}{\sqrt{2}}(e_1\wedge e_2-e_3\wedge e_4)$ have the following form in the basis (\ref{X1})
\begin{eqnarray*}
g_1=\beta_1 X_1+\beta_2 X_2+\beta_3 X_3, \ \ \
g_2=\gamma_1Y_1+\gamma_2 Y_2+\gamma_3 Y_3,
\end{eqnarray*}
where  
\begin{eqnarray*}
\beta_1= \frac{1+a_{10}b_{01}-a_{01}b_{10}}{\sqrt{2}\ h}, \ \beta_2=\frac{b_{01}-a_{10}}{\sqrt{2}\ h}, \ \beta_3=-\frac{a_{01}+b_{10}}{\sqrt{2}\ h}, \label{beta}\\
\gamma_1= \frac{1-a_{10}b_{01}+a_{01}b_{10}}{\sqrt{2}\ h}, \ \gamma_2=-\frac{a_{10}+b_{01}}{\sqrt{2}\ h}, \ \gamma_3=\frac{b_{10}-a_{01}}{\sqrt{2}\ h}.\label{gamma}
\end{eqnarray*} 

It is easy to see that these functions satisfy (\ref{2-sphere}). By applying the stereographic projection to $(\beta_{1},\beta_2,\beta_3)\in S_1$,
\begin{equation*}
(\beta_{1},\beta_2,\beta_3)\mapsto \left(\frac{\beta_2}{\frac{1}{\sqrt{2}}+\beta_1},\frac{\beta_3}{\frac{1}{\sqrt{2}}+\beta_1}\right)\in \mathbb R^2
\end{equation*}
and analogously to $(\gamma_1,\gamma_2,\gamma_3)\in S_2$, we get that the components of the Gauss map have the following local form

\begin{equation}\label{g1g2}
  \begin{array}{ccl}
g_1(u,v)=r\left(-a_{10}+b_{01},-a_{01}-b_{10}\right), \ \ \  
	g_2(u,v)=p\left(-a_{10}-b_{01},b_{10}-a_{01}\right), 
  \end{array}
\end{equation}
where
\begin{equation}\label{rr}
r=\frac{1}{ 1-a_{01}b_{10}+a_{10}b_{01}+h} \ \ p=\frac{1}{ 1+a_{01}b_{10}-a_{10}b_{01}+h},
\end{equation}
and $h$ is as (\ref{h}).
We notice that for any $(u,v)$
\begin{eqnarray*} \label{positivity1}
1-a_{01}b_{10}+a_{10}b_{01}+h>1,\ \ \
1+a_{01}b_{10}-a_{10}b_{01}+h>1,\nonumber 
\end{eqnarray*}
 since 
\begin{equation*}
h> \left|a_{01}b_{10}-a_{10}b_{01}\right|.
\end{equation*}
These inequalities imply that the components of the Gauss map in (\ref{g1g2}) are smooth map-germs  $\mathbb R^2\rightarrow \mathbb R^2$.
To simplify notation, we denote the components of these map-germs by $g_1$ and $g_2$.   

\section{Generic singularities of the components of the Gauss map} \label{sing}

In \cite{L} Little proved that the Gauss map of an immersion $M$ into $\mathbb R^4$ is not an immersion at $p$ if and only if $\Delta=0$ and $K=0$ at $p$.   
We begin by showing that the Gauss map is generically an immersion in the setting. Namely,
\begin{prop}\label{codim3}
If $M$ is a generic,
closed,
oriented surface immersed into $\mathbb R^4$ then the Gauss map $g:M\rightarrow  Gr^{+}(2,4)$ is an immersion.
\end{prop}

\begin{proof}
For simplicity, we assume that an immersed surface $M$ into $\mathbb R^4$ is parameterized locally in Monge form at the origin. Namely,
it is
locally described as a graph of a smooth map-germ 
\begin{eqnarray}\label{graf}
F:\mathbb R^2 \rightarrow \mathbb R^2,\ \ \ (u,v)\mapsto (a(u,v),b(u,v)), 
\end{eqnarray}
such that
\begin{equation*}\label{Monge_at_0}
a_{10}(0,0)=a_{01}(0,0)=b_{10}(0,0)=b_{01}(0,0)=0.
\end{equation*}
Thus, we use the Gauss map's local expression (\ref{g1g2}) of $$g=(g_1,g_2):\mathbb R^2 \rightarrow \mathbb R^2\times \mathbb R^2$$ to 
develop the analysis of the space of jets.

Let $J^k(2,2)$ be the space of $k$-jets of smooth maps $C^{\infty}(\mathbb R^2,\mathbb R^2)$.  
The Gauss map is not an immersion at $(0,0)$ if its differential 
\begin{equation}\label{mdg}
dg_{(0,0)}=\left(\begin{array}{cc}
S^1_{10}(0,0) & S^1_{01} (0,0)\\
T^1_{10}(0,0) & T^1_{01} (0,0)\\
S^2_{10}(0,0) & S^2_{01}(0,0) \\
T^2_{10} (0,0)& T^2_{01}(0,0) \\
\end{array}
\right)
\end{equation}
has a rank $0$ or $1$, where
\begin{eqnarray*}
S^1_{10}=r_{00}(a_{20}-b_{11})+r_{10}(a_{10}-b_{01}),&&
S^1_{01}=r_{00}(a_{11}-b_{02})+r_{01}(a_{10}-b_{01}),\\
T^1_{10}=r_{00}(a_{11}+b_{20})+r_{10}(a_{01}+b_{10}),&&
T^1_{01}=r_{00}(a_{02}+b_{11})+r_{01}(a_{01}+b_{10}), \\
S^2_{10}=p_{00}(a_{20}+b_{11})+p_{10}(a_{10}+b_{01}),&&
S^2_{01}=p_{00}(a_{11}+b_{02})+p_{01}(a_{10}+b_{01}),\\
T^2_{10}=p_{00}(a_{11}-b_{20})+p_{10}(a_{01}-b_{10}),&&
T^2_{01}=p_{00}(a_{02}-b_{11})+p_{01}(a_{01}-b_{10}),
\end{eqnarray*}

\begin{eqnarray}
&h_{00}=\sqrt{1+a_{10}^2+a_{01}^2+b_{10}^2+b_{01}^2+(a_{01}b_{10}-a_{10}b_{01})^2},&\nonumber\\
&r_{00}=1/\left( 1-a_{01}b_{10}+a_{10}b_{01}+h_{00}  \right),& \label{R00}\\
&r_{10}=-r_{00}^2(a_{20} b_{01}-a_{11} b_{10}+a_{10} b_{11}-a_{01} b_{20}+&\nonumber \label{R10}\\
&\frac{a_{01} a_{11}+a_{10}
	a_{20}+b_{01} b_{11}+b_{10} b_{20}+(a_{10} b_{01}-a_{01} b_{10}) (a_{20} b_{01}-a_{11} b_{10}+a_{10}
	b_{11}-a_{01} b_{20})}{h_{00}}),& \nonumber \\
&r_{01}=-r_{00}^2(a_{11} b_{01}+a_{10} b_{02}-a_{02} b_{10}-a_{01} b_{11}+&\nonumber \label{R01}\\
&\frac{a_{01} a_{02}+a_{10} a_{11}+b_{01}
	b_{02}+b_{10} b_{11}+(a_{10} b_{01}-a_{01} b_{10}) (a_{11} b_{01}+a_{10} b_{02}-a_{02} b_{10}-a_{01}
	b_{11})}{h_{00}}),& \nonumber \\
&p_{00}=1/\left( 1+a_{01}b_{10}-a_{10}b_{01}+h_{00}  \right),& \\
&p_{10}=-p_{00}^2(-a_{20} b_{01} + a_{11} b_{10} - a_{10} b_{11} + a_{01} b_{20} +&\nonumber \label{P10}\\
 &\frac{a_{01} a_{11} + a_{10} a_{20} + b_{01} b_{11} + 
      b_{10} b_{20} + (a_{10} b_{01} - a_{01} b_{10}) (a_{20} b_{01} - a_{11} b_{10} + a_{10} b_{11} - 
         a_{01} b_{20})}{h_{00}}),& \nonumber \\
&p_{01}=-p_{00}^2(-a_{11} b_{01} - a_{10} b_{02} + a_{02} b_{10} + a_{01} b_{11} +&\nonumber \label{P_{01}}\\
 &\frac{a_{01} a_{02} + a_{10} a_{11} + b_{01} b_{02} + 
     b_{10} b_{11} + (a_{10} b_{01} - a_{01} b_{10}) (a_{11} b_{01} + a_{10} b_{02} - a_{02} b_{10} - 
        a_{01} b_{11})}{h_{00}}).& \nonumber
\end{eqnarray}

If the rank of $dg_{(0,0)}$ is $1$, then at least one of the rows of  matrix (\ref{mdg}) does not vanish. For simplicity, we  assume that 
$(S^1_{10},S^1_{01})\ne (0,0)$.
Then, the rank of ${dg}_{(0,0)}$ is $1$ if and only if 
$$
W_1=S^1_{10}T^1_{01}-S^1_{01}T^1_{10}, \ W_2=S^1_{10}S^2_{01}-S^1_{01}S^2_{10}, \ W_3=S^1_{10}T^2_{01}-S^1_{01}T^2_{10}
$$
vanish at $(0,0)$.
We consider the following subset 
$$\Gamma^1=\{(S^1_{10})^2+(S^1_{01})^2\ne 0, W_1=W_2=W_3=0\}$$
of $J^3(2,2)$.
A direct calculation implies that 
$$dW_1\wedge dW_2 \wedge dW_3\ne 0\  {\rm for}\ (S^1_{10})^2+(S^1_{01})^2\ne 0
\ {\rm and}\ a_{10}=a_{01}=b_{10}=b_{01}=0.$$ 
Thus, $\Gamma^1$ is a submanifold of codimension $3$. So, for a generic $F$, $j^3F$ does not intersect this set.

If the rank of  ${dg}_{(0,0)}$ is $0$, then 
\begin{equation}\label{rank-zero}
S^i_{10}=T^i_{10}=S^i_{01}=T^i_{01}=0, \ \text{for \  i=1,2\ at\ }(0,0). 
\end{equation}
 It is easy to check that a subset of $J^3(2,2)$ described by (\ref{rank-zero}) is contained in the following set 
 \begin{equation*}\label{rank-zero1}
C=\left\lbrace S^1_{10}=T^1_{10}=S^2_{10}=T^2_{10}=S^1_{01}=T^1_{01}=0\right\rbrace . 
\end{equation*}
Moreover, $C$ is a submanifold of codimension $6$ in $J^2(2,2)$, because 
$$
dS^1_{10}\wedge dT^1_{10}\wedge dS^2_{10}\wedge dT^2_{10}\wedge dS^1_{01}\wedge dT^1_{01}\ne 0
$$
for $a_{10}=a_{01}=b_{10}=b_{01}=0$.
Hence, for a generic $F$, $j^2F$ does not intersect this set. Finally, we point out that this local argument extends in the standard way for immersions of closed surfaces into $\mathbb R^4$.
\end{proof}

Now we state the main theorem.
\begin{thm}\label{generic}
Properties $(G_1)$, $(G_2)$ and $(G_3)$ are generic.
\end{thm} 
\begin{proof} 
We only present a proof for the first component $g_1$ of the Gauss map since a proof for the second component $g_2$ is similar.
We use the local descriptions  (\ref{g1g2}) and (\ref{graf}) and employ Wolfram Mathematica \cite{Math} for the more complicated computations.
Let $J^k(2,2)$ be the space of $k$-jets of smooth maps $C^{\infty}(\mathbb R^2,\mathbb R^2)$.  
The 4-jet of $F$ at $(u,v)$ has the following form
 \begin{eqnarray*}
j^4_{(u,v)} F&=&\left(\begin{array}{c}\sum_{i+j=0}^4\frac{a_{ij}}{i!j!}\Delta u^i \Delta v^j\\
\sum_{i+j=0}^4\frac{b_{ij}}{i!j!}\Delta u^i \Delta v^j
\end{array}\right),
\end{eqnarray*}
where $\Delta u=(\tilde{u}-u)$ and $\Delta v=(\tilde{v}-v)$.
 By (\ref{g1g2}),  the  3-jet of $g_1$ at $(u,v)$ has the following form
\begin{eqnarray}
j^3_{(u,v)}g_1&=&\left(\begin{array}{c}\sum_{i+j=0}^3\frac{s_{ij}}{i!j!}\Delta u^i \Delta v^j\\
\sum_{i+j=0}^3\frac{t_{ij}}{i!j!}\Delta u^i \Delta v^j
\end{array}\right), \label{j^3g1}
\end{eqnarray}
where
\begin{equation}\label{s}
  \begin{array}{ccl}
s_{kl}&=&-\frac{\partial^{k+l}}{\partial u^k \partial v^l}\left (r(u,v)(a_{10}(u,v)-b_{01}(u,v))\right), \\
 t_{kl}&=&-\frac{\partial^{k+l}}{\partial u^k \partial v^l} (r(u,v)(a_{01}(u,v)+b_{10}(u,v))).
  \end{array}
\end{equation}

We apply formulas (\ref{j^3g1})-(\ref{s}) to define the map 
$$\mu:J^4(2,2)\rightarrow J^3(2,2), \ \ j^4_{(u,v)}F\mapsto j^3_{(u,v)}g_1.$$ 
Then, we will show that  $\mu$ is a submersion. 
In the natural coordinate system  
$$
x=\left(u,v,a_{ij},b_{ij}\right)_{i+j=0,\cdots,4}, \ y=\left(u,v,s_{ij},t_{ij}\right)_{i+j=0,\cdots,3}$$
on $J^4(2,2)$ and $J^3(2,2)$, respectively, the map $\mu$ has the following form
\begin{equation*}\label{gamma1}
\mu(x)=(u,v,S_{00}(x),T_{00}(x),S_{10}(x),T_{10}(x),\cdots,S_{03}(x),T_{03}(x)),
\end{equation*}
where by (\ref{g1g2}),      
we have
\begin{eqnarray*} 
S_{00}(x)=-r_{00}(a_{10}-b_{01}),\label{S00}\\
T_{00}(x)=-r_{00}(a_{01}+b_{10}),\label{T00}
\end{eqnarray*}
$r_{00}$ is given by (\ref{R00}), 
and $S_{ij}(x), \ T_{ij}(x)$ are obtained by applying to formulas (\ref{s}) the following rules 
\begin{equation} \label{rules}
\frac{\partial^{k+l}}{\partial u^k \partial v^l} a_{ij}=a_{(i+k)(j+l)}, \ \frac{\partial^{k+l}}{\partial u^k \partial v^l} b_{ij}=b_{(i+k)(j+l)} .
\end{equation}
 For example, the functions $S_{10}, T_{10}, S_{01}, T_{01}$ have the following form
{\small
\begin{equation*}
\label{dg1}\begin{array}{cc}	
S_{10}(x)=-r_{00}(a_{20}-b_{11})-r_{10}(a_{10}-b_{01}),\nonumber& S_{01}(x)=-r_{00}(a_{11}-b_{02})-r_{01}(a_{10}-b_{01}),\\
T_{10}(x)=-r_{00}(a_{11}+b_{20})-r_{10}(a_{01}+b_{10}),& T_{01}(x)=-r_{00}(a_{02}+b_{11})-r_{01}(a_{01}+b_{10}),
\end{array}
\end{equation*}}
where $r_{10}$ and $r_{01}$ are given by (\ref{R00}). Observe that these formulas coincide with those in $(\ref{mdg} )$,
after omitting superindex $1$. 
The explicit formulas for $S_{ij}$, $T_{ij}$ are relatively easy to calculate, but they should be shorter to present all of them.  
For a smooth map $(w_1,\cdots,w_n)\mapsto (P_1(w_1,\cdots,w_n),\cdots,P_m(w_1,\cdots,w_n))$, we denote by 
$\frac{\partial(P_{i_1},\cdots,P_{i_k})}{\partial
(w_{j_1},\cdots,w_{j_k})}$ the matrix $\left(\frac{\partial P_{i_l}}{\partial w_{j_q}}\right)_{l,q=1,\cdots,k}$.
Let us notice that $T_{ij}$ are $S_{ij}$ do not depend on $a_{kl}$ and $b_{kl}$ for $k+l>i+j+1$. Hence, the square
 matrix of first order partial derivatives of $\mu$ concerning variables $u,v,a_{ij},b_{ij}$ for $i\ne 0$ has the following form
\begin{eqnarray}
&\frac{\partial \mu}{\partial(u,v,a_{ij},b_{ij})}_{i\ne 0}=&\label{matrix}\\ 
&\left(\begin{array}{ccccc} \nonumber
I_2&0&0&0&0\\
\ast&\frac{\partial(S_{00},T_{00})}{\partial(a_{10},b_{10})}&0&0&0\\
\ast&\ast&\frac{\partial(S_{10},T_{10},S_{01},T_{01})}{\partial(a_{20},b_{20},a_{11},b_{11})}&0&0\\
\ast&\ast&\ast&\frac{\partial(S_{20},T_{20},\cdots,S_{02},T_{02})}{\partial(a_{30},b_{30},\cdots,a_{12},b_{12})}&0\\
\ast&\ast&\ast&\ast&\frac{\partial(S_{30},T_{30},\cdots,S_{03},T_{03})}{\partial(a_{40},b_{40},\cdots,a_{13},b_{13})}
\end{array}\right).
\end{eqnarray}
By direct computations, we get
$$
\det \frac{\partial(S_{k0},T_{k0},\cdots,S_{0k},T_{0k})}{\partial(a_{(k+1)0},b_{(k+1)0},\cdots,a_{1k},b_{1k})}=\epsilon_k \left(\frac{1+a_{01}^2+b_{10}^2}{h_{00}}r_{00}^2\right)^{k+1},$$ 
where $\epsilon_k=-1$ for $k=1,2$ and $\epsilon_k=1$ for $k=0,3$. Thus by (\ref{matrix}), we obtain that 
$$\det \frac{\partial \mu}{\partial(u,v,a_{ij},b_{ij})}_{i\ne 0}=\left(\frac{1+a_{01}^2+b_{10}^2}{h_{00}}r_{00}^2\right)^{10}>0.$$ Hence,  $\mu:J^4(2,2)\rightarrow J^3(2,2)$ is a submersion. 
It implies that for any smooth submanifold $\Sigma$ of $J^3(2,2)$, the set $\mu^{-1}(\Sigma)$ is a smooth 
submanifold of $J^4(2,2)$ and $\text{codim}\ \mu^{-1}(\Sigma)=\text{codim}\ \Sigma$.  The map 
$j^3g_1:
(u,v)\mapsto j^3_{(u,v)}g_1\in J^3(2,2)$ is transverse to $\Sigma$ if and only if the map 
$j^4F:
(u,v)\mapsto j^4_{(u,v)}F\in J^4(2,2)$ is tranverse to $\mu^{-1}(\Sigma)$ of $J^4(2,2)$ 
since $j^3g_1=\mu \circ j^4F$ and $\mu$ is a submersion.

Let $dg_1=\left(\begin{array}{c}
s_{10}, \ s_{01} \\ 
t_{10}, \ t_{01} 
\end{array}\right)$ 
and 
$$
Q_{g_1}(\Delta u, \Delta v)=\left(\begin{array}{c}
s_{20}\Delta u^2+ 2s_{11}\Delta u\Delta v+ s_{02}\Delta v^2 \\ 
t_{20}\Delta u^2+ 2t_{11}\Delta u\Delta v+ t_{02}\Delta v^2 
\end{array}\right).
$$
In this notation, the  3-jet of $g_1$ at $(u,v)$ has the following form
\begin{eqnarray*}
j^3_{(u,v)}g_1&=&
\left(\begin{array}{c} s_{00}\\ 
t_{00} 
\end{array}\right)+
dg_1
\left(\begin{array}{c} 
\Delta u\\
\Delta v
\end{array}\right)+\frac{1}{2}Q_{g_1}(\Delta u, \Delta v)\nonumber \\
&&+\frac{1}{6}\left(\begin{array}{c}
s_{30}\Delta u^3+ 3s_{21}\Delta u^2\Delta v+ 3s_{12}\Delta u\Delta v^2 +s_{03}\Delta v^3\\ 
t_{30}\Delta u^3+ 3t_{21}\Delta u^2\Delta v+ 3t_{12}\Delta u\Delta v^2 +t_{03}\Delta v^3
\end{array}\right).
\end{eqnarray*}
We consider $\Sigma^i=\{\dim\ker dg_1=i\}\subset J^4(2,2)$. Then, $J^4(2,2)=\bigcup_{i=0}^2 \Sigma^i$ and  $\Sigma^i$ is a smooth submanifold of codimension $i^2$. So, the point $(u,v)$ is a regular point of $g_1$ if and only if $j^3_{(u,v)}g_1$ belongs to $\Sigma^0=\{s_{10}t_{01}-s_{01}t_{10}\ne 0\}$. 

The submanifold $\Sigma^1$ can be described in the following way 
\begin{equation*}\label{Sigma1}
\Sigma^1=\{s_{10}t_{01}-s_{01}t_{10}=0, \ s_{10}^2+t_{10}^2+s_{01}^2+t_{01}^2\ne 0\}.
\end{equation*} 
We observe that the immersion $(u,v)\mapsto (u,v,F(u,v))$ has property $(G_1)$ if and only if $j^3g_1$ is transverse to $\Sigma^1$. Then the set 
$\{(u,v)| j^3_{(u,v)}g_1\in \Sigma^1\}$ is a smooth submanifold of codimension $1$. 
The point $(u,v)$ is a singular point of $g_1$ of corank $1$ if and only if 
$j^3_{(u,v)}g_1$ belongs to  $\Sigma^1$. The singular point $(u,v)$ of corank $1$ is a fold point if the kernel of ${dg_1}_{(u,v)}$ is transverse to the submanifold 
$\{(u,v)| j^3_{(u,v)}g_1\in \Sigma^1\}$. It means that the restriction of $g_1$ to the smooth curve of singular points is a regular map 
at $(u,v)$; or equivalently (see Section 7.1 of \cite{Martinet})
\begin{equation}\label{fold}
Q_{g_1}(\ker {dg_1}_{(u,v)})\not\subset {dg_1}_{(u,v)}(\mathbb R^2).
\end{equation}
Since  $s_{10}^2+t_{10}^2+s_{01}^2+t_{01}^2\ne 0$, we may assume that $s_{10}\ne 0$. In other cases the proofs are similar. Since  $s_{10}\ne 0$ and $s_{10}t_{01}-s_{01}t_{10}=0$, $\ker dg_1$ is spanned by the vector $\left(-s_{01},s_{10}\right)$, and the image of $dg_1$ is spanned by $\left(s_{10},t_{10}\right)$. Then, condition (\ref{fold}) has the following form
\begin{equation*}
\det\left(\begin{array}{cc} s_{20}s_{01}^2- 2s_{11}s_{01}s_{10}+ s_{02}s_{10}^2, & s_{10} \\ 
t_{20}s_{01}^2- 2t_{11}s_{01}s_{10}+ t_{02}s_{10}^2, & t_{10}
\end{array}          \right)\ne 0.
\end{equation*}
We denote by $\Sigma^{11}\subset \Sigma^1$ the following subset of $J^3(2,2)$  
\begin{equation*}\label{Sigma11}
\Sigma^{11}=\{ \dim\ker dg_1=1, \ Q_{g_1}(\ker dg_1)\subset dg_1(\mathbb R^2)\}.
\end{equation*} 
 Since  $s_{10}\ne 0$,   
\begin{equation*}\Sigma^{11}=\left\{s_{10}t_{01}-s_{01}t_{10}=0, \ 
\det\left(\begin{array}{cc} s_{20}s_{01}^2- 2s_{11}s_{01}s_{10}+ s_{02}s_{10}^2, & s_{10} \\ 
t_{20}s_{01}^2- 2t_{11}s_{01}s_{10}+ t_{02}s_{10}^2, & t_{10}
\end{array}          \right)= 0
\right\} 
\end{equation*}
is a smooth submanifold of codimension $2$ of $J^4(2,2)$. Therefore, the immersion
 $(u,v)\mapsto (u,v,F(u,v))$ has property $(G_2)$ if and only if $j^3g_1$ is transverse to $\Sigma^1$ and $\Sigma^{11}$.

Let us denote $\Sigma^{10}=\Sigma^{1}\setminus \Sigma^{11}$. The point $(u,v)$ is a fold point of $g_1$ if and only if $j^3_{(u,v)}g_1$ belongs to $\Sigma^{10}$ (see Section 7.3 of \cite{Martinet}). The point $(u,v)$ is a cusp point of $g_1$ if and only if $j^3_{(u,v)}g_1$ belongs to $\Sigma^{11}$ and the map $j^3g_1$ is transverse to $\Sigma^{11}$ at $(u,v)$.  By Thom's Transversality Theorem, the map $g_1$ is generic if $j^3g_1$ is transverse to submanifolds $\Sigma^0$, $\Sigma^{10}$, $\Sigma^{11}$, $\Sigma^{2}$.  Since the codimension of $\Sigma^2$ is $4$, then $j^3g_1$ is transverse to $\Sigma^2$ if its image does not intersect $\Sigma^2$. Hence, generically singularities of $g_1$ are fold points and cusp points (see Theorem 7.4 
in \cite{Martinet}). 

Since $\mu$ is a submersion,  $j^3g_1$ is transverse to the submanifold $\Sigma^0$ ($\Sigma^{10}$, $\Sigma^{11}$, $\Sigma^{2}$ respectively)  if and only if $j^4F$ is transverse to the submanifold $\mu^{-1}(\Sigma^0)$ ($\mu^{-1}(\Sigma^{10})$, 
$\mu^{-1}(\Sigma^{11})$, $\mu^{-1}(\Sigma^{2})$ respectively). The set of points of $M$ at which $j^4F$ transversally intersects $\mu^{-1}(\Sigma^1)$ is a $1$-dimensional smooth submanifold of the surface $M$.  Since $M$ is closed, it is a finite disjoint collection of  smoothly embedded circles. Since $\text{codim}\ \mu^{-1}(\Sigma^{11})=2$, the set of points of $M$ at which $j^4F$  intersects transversally  $\mu^{-1}(\Sigma^{11})$ is a $0$-dimensional smooth submanifold of the surface $M$. So, it
consists of a finite set of points (cusp points) on the singular set of $g_1$.  All the other singular points of $g_1$ are fold points because $\text{codim}\ \mu^{-1}( \Sigma^2)=4$.  
Thus, by Thom's Transversality Theorem, we finish the proof of genericity of $(G_1)$ and $(G_2)$.

To prove the genericity of $(G_3)$, we first define the sets 
$$
M^{(r)}=\{(p_1,\cdots, p_r)\in M^r|p_i\ne p_j \ \text{for} \ i\ne j\},
$$ 
$$
\Delta_r=\{(q,\cdots,q)\in N^r|q\in N\},
$$
where $M$ and $N$ are smooth manifolds. Let
$\alpha:J^k(M,N)\rightarrow M$ be the canonical projection. Then, we define
$$
\alpha^r=\alpha\times \cdots \times \alpha:(J^k(M,N))^r\rightarrow M^r,\ \ (j^k_{p_1} f_1,\cdots,j^k_{p_r} f_r)\mapsto (p_1,\cdots, p_r),
$$
and determine a $r$-fold $k$-jet bundle as 
$_rJ^k(M,N)=(\alpha^r)^{-1}(M^{(r)})$. Thus, for a map $f\in C^{\infty}(M,N)$, we define 
$$
_rj^kf:M^{(r)}\rightarrow _rJ^k(M,N),\ \ (p_1,\cdots,p_r)\mapsto (j^k_{p_1}f,\cdots,j^k_{p_r}f).
$$
The $P=(\mathbb R^2)^{(2)}\times \Delta_2\times (\Sigma^1)^2$ is a submanifold of $_2J^3(2,2)$ of codimension $4$, $C=(\mathbb R^2)^{(2)}\times \Delta_2\times \Sigma^1\times \Sigma^{1,1}$ is a submanifold of $_2J^3(2,2)$ of codimension $5$, and  $T=(\mathbb R^2)^{(3)}\times \Delta_3\times (\Sigma^1)^3$ is a submanifold of $_3J^3(2,2)$ of codimension $7$. Thus, property $(G_3)$ means that $_2j^3g_1$ is transverse to $P$, and $C$ and $_3j^3g_1$ is transverse to $T$. Since $\text{codim} \ C$ is $5$ and $\dim (\mathbb R^2)^{(2)}$ is $4$, the image of  $_2j^3g_1$ does not intersect $C$. Moreover, since $\text{codim} \ T$ is $7$ and $\dim (\mathbb R^2)^{(3)}$ is $6$,  the image of  $_3j^3g_1$ does not intersect $T$. The map $_2j^3g_1$ is the restriction to $(\mathbb R^2)^{(2)}$ of the map 
$$
j^3g_1\times j^3g_1:\mathbb R^2\times \mathbb R^2\rightarrow J^3(2,2)\times J^3(2,2).
$$
Since $P$ does not intersect $(j^3g_1\times j^3g_1)(\Delta_2)$, the map $_2j^3 g_1$ is transverse to $P$ if and only if $j^3g_1\times j^3g_1$ is transverse to $P$, and so 
$(_2j^3g_1)^{-1}(P)=(j^3g_1\times j^3g_1)^{-1}(P)$ is a finite collection of pairs of points since $M \times M$ is compact. 

 Let $\pi_i:(\mathbb R^2)^{(2)}X \rightarrow \mathbb R^2,\ \ (p_1,p_2) \mapsto p_i$ for $i=1,2$. Then, the transverse intersection of 
$_2j^3g_1$ and  $P$  is equivalent to the non-vanishing of the following $4$-form 
\begin{equation}\label{4-form}
d(\pi_1^{\ast}g_{11}-\pi_2^{\ast}g_{11})\wedge d(\pi_1^{\ast}g_{12}-\pi_2^{\ast}g_{12})\wedge
 d(\pi_1^{\ast}(J(g_1)))\wedge d(\pi_2^{\ast}(J(g_1)))
\end{equation}
 at $(q_1,q_2)\in (\mathbb R^2)^{(2)}$ such that $g_1(q_1)=g_1(q_2)$,
 where $g_1=(g_{11},g_{12})$. Let $X_i$ be a vector spanning the kernel of ${dg_1}_{q_i}$ and $Y_i$ be a vector spanning the tangent space to $J(g_1)$ at the point $q_i$ for $i=1,2$. Since $q_1\ne q_2$ and $q_1$, $q_2$ are fold points, vectors $X_1,Y_1\in T_{q_1}\mathbb R^2$, $X_2,Y_2\in T_{q_2}\mathbb R^2$ span the tangent space $T_{(q_1,q_2)}(\mathbb R^2)^{(2)}=T_{q_1}\mathbb R^2\oplus T_{q_2}\mathbb R^2$;  thus, for a smooth function $f$ on $\mathbb R^2$ and $i\ne j$, we have 

\begin{equation}\label{equal}
  \begin{array}{ccl}
 {d(\pi_i^{\ast}f)}_{(q_1,q_2)}(X_j)={d(\pi_i^{\ast}f)}_{(q_1,q_2)}(Y_j)=0,\\  
{d(\pi_i^{\ast}f)}_{(q_1,q_2)}(X_i)={df}_{q_i}(X_i),\ \ {d(\pi_i^{\ast}f)}_{(q_1,q_2)}(Y_i)={df}_{q_i}(Y_i). 
  \end{array}
\end{equation}
By evaluating the $4$-form (\ref{4-form}) on $(X_1,X_2,Y_1,Y_2)$, we obtain by (\ref{equal}) that
 $$
 \det[ {dg_{1}}_{q_2}(Y_2),{dg_{1}}_{q_1}(Y_1)]{d(J(g_1))}_{q_1}(X_1){d(J(g_1))}_{q_2}(X_2)\ne 0,
 $$
 which implies that vectors ${dg_{1}}_{q_1}(Y_1)$ and ${dg_{1}}_{q_2}(Y_2)$ are linearly independent.  
 Thus, the image of the singular set intersects transversally at  $g_1(q_1)=g_1(q_2)$.
 
 For a submanifold $Q$ of $_rJ^3(2,2)$, $_rj^4 F$ is transverse to $(\mu^r)^{-1}(Q)$ if and only if $_rj^3 g_1$ is transverse to $Q$. Thus, by the Multijet Transversality Theorem (see Theorem 4.13 in Chapter II \cite{GG}), property $(G_3)$ is generic.  
\end{proof}

By Theorem 6.3 in Chapter VII \cite{GG} (see also Proposition 2.4 in \cite{BW}), we get the following result.
\begin{cor}
Let $M$ be an 
oriented, closed surface.
Let $M\rightarrow \mathbb R^{4}$ be an embedding that has property $(G_3)$. 
If $g: M\rightarrow S_{1}\times S_{2}$ is the Gauss map  
then, the component $g_i,\ i=1,2$ is stable.
 \end{cor} 



\section{Singularities of $g_i$ and contact with holomorphic curves.}

Singularities of the Gauss map of a surface in $\mathbb R^3$ are studied by contact of the surface with its tangent planes. The Gauss map of an affine plane is constant, so its tangent map vanishes at every point.
Analogously, we study the singularities of the Gauss map components of a surface in $\mathbb R^4$ by the contact of this surface with $\mathcal J$-holomorphic curves. To do so, we first characterize a $\mathcal J$-holomorphic curve, with respect to an orthogonal complex structure $\mathcal J$ on $\mathbb{R}^4$, as a surface with a constant Gauss map component (see Theorem 5.3 (b) in \cite{HO1}).

We recall that an orthogonal complex structure on the euclidean space $\mathbb R^4$ is an integrable almost complex structure $\mathcal J$ on $\mathbb R^4$ such that $<\mathcal Jv,\mathcal Jw>=<v,w>$ for any $x\in \mathbb R^4$ and any $v,w\in T_x \mathbb R^4$. 
An example of an orthogonal complex structure on $\mathbb R^4$ is the standard complex structure on $\mathbb R^4=\mathbb C^2$ denoted by $\mathcal J_0$. 
It is worth noticing that every orthogonal complex structure on $\mathbb R^4$  is constant i.e. $d \mathcal J=0$ (see proof of Proposition 6.6 in \cite{W} and Corollary 3.5 in \cite{Ch}).  

\begin{prop}
A connected surface $M$ is a $\mathcal J$-holomorphic curve with respect to an orthogonal complex structure $\mathcal J$ on an euclidean space $\mathbb R^4$ if and only if 
$g_i,\ i=1\ {\rm or}\ 2$ is constant. If  i=1, $\mathcal J$ preserves the orientation of  $\mathbb R^4$; otherwise, it reverses it. Moreover, in any case, the mean curvature $H$ vector vanishes at any point. 
\end{prop}

\begin{proof} Let us assume that the orthogonal complex structure $\mathcal J$ preserves the orientation of $\mathbb R^4$. Since $M$ is a $\mathcal J$-holomorphic curve, we can choose a Darboux frame $(e_1, e_2, e_3, e_4)$ with $e_2=\mathcal Je_1$ and $e_4= \mathcal J e_3$. Since $\mathcal J$ is constant, we have
$$
\omega_2^3=<de_2,e_3>=<d(\mathcal J e_1),e_3>=<\mathcal Jde_1,e_3>=
$$
$$
<\mathcal Jde_1,-\mathcal J^2e_3>=-<de_1,\mathcal Je_3>=-<de_1,e_4>=-\omega_1^4,
$$
and
$$
\omega_2^4=<de_2,e_4>=<d(\mathcal Je_1),\mathcal Je_3>=<\mathcal Jde_1,\mathcal Je_3>=<de_1,e_3>=\omega_1^3.
$$
Thus, we obtain that 
\begin{equation}\label{dg_1=0}
\omega_2^3+\omega_1^4=\omega_2^4-\omega_1^3=0.
\end{equation}
It implies that $dg_1=0$, by (\ref{dg_1}).
Conversely, if $g_1$ is constant, we use the local expression of $g_1$ in $(\ref{g1g2})$ to conclude, after evaluating $(u,v)=(0,0)$, that the map $(a(u,v),b(u,v))$
satisfies the Cauchy-Riemann conditions. Now, we observe that
the mean curvature vector $H$ is given by 
\begin{equation}\label{mean}
H=\frac{1}{2}\left((\omega_1^3(e_1)+\omega_2^3(e_2))e_3+(\omega_1^4(e_1)+\omega_2^4(e_2))e_4\right).
\end{equation}
 By (\ref{dg_1=0}) and (\ref{w1j2=w2j1}), we get that
$$
\omega_1^3(e_1)=\omega_2^4(e_1)=\omega_1^4(e_2)=-\omega_2^3(e_2)
$$
and 
$$
\omega_1^4(e_1)=-\omega_2^3(e_1)=-\omega_1^3(e_2)=-\omega_2^4(e_2).
$$
Thus, we get that $\omega_1^3(e_1)+\omega_2^3(e_2)=\omega_1^4(e_1)+\omega_2^4(e_2)=0$, which implies that $H=0$ by (\ref{mean}). 
Suppose the orthogonal complex structure reverses the orientation of $\mathbb R^4$ and $M$ is a holomorphic curve. In that case, we can choose a Darboux frame $(e_1, e_2, e_3, e_4)$  such that $e_2=\mathcal Je_1$ and $e_4=-\mathcal Je_3$ and proceed in the same way described above.

\end{proof} 

We recall some basic definitions and results on contact equivalence (for details
see \cite{AGV}, \cite{DRR}).

\begin{defn}\label{Keq}
The map-germs $f, \widetilde{f}:(\mathbb R^s,0)\rightarrow (\mathbb
R^t,0)$ are $\mathcal K$-{\bf equivalent}  if there
exists a diffeomorphism-germ $\phi:(\mathbb R^s,0)\rightarrow
(\mathbb R^s,0)$ and a map-germ $A:(\mathbb R^s,0)\rightarrow
GL(\mathbb R^t)$ \ such that \
$
\tilde f=A\cdot (f\circ \phi).
$
\end{defn}

Let $\mathcal E_s$ denote the local ring of smooth function-germs at $0$ on $\mathbb R^s$.

\begin{rem}[\cite{AGV}] \label{V-ideal}
For the $\mathcal K$-equivalence of two map-germs it is necessary
and sufficient that the two ideals generated by the components of
these map-germs may be mapped one to the other by an
isomorphism of $\mathcal E_s$ induced by a  diffeomorphism-germ of the source space $(\mathbb R^s,0)$.
\end{rem}
Let $M_i, N_i$ be germs at $x_i$ for $i=1,2$ of smooth surfaces in
the space $\mathbb R^{4}$. We describe them in the following way.

\noindent $M_i=\phi_i^{-1}(0)$, where $\phi_i:(\mathbb R^4,x_i)\rightarrow (\mathbb R^{2},0)$ is a submersion-germ for $i=1,2$, and \\
\noindent $N_i=\psi_i(\mathbb R^2)$, where $\psi_i:(\mathbb R^2,0)\rightarrow (\mathbb R^4,x_i)$ is an embedding-germ for $i=1,2$.

\begin{defn}
The contact of $M_1$ and $N_1$ at $x_1$ is of the same {\bf contact-type} as the contact of $M_2$ and $N_2$ at $x_2$ if there exists a diffeomorphism-germ
$\Phi:(\mathbb R^4,x_1)\rightarrow (\mathbb R^4,x_2)$ such that $\Phi(M_1)=M_2$ and $\Phi(N_1)=N_2$. We denote the contact-type of $M_1$ and $N_1$ at $x_1$ by $\mathcal K(M_1,N_1,x_1)$.
\end{defn}

The contact-type can be studied by a contact-map.
 
\begin{defn}
A {\bf contact map} between surface-germs $M_1, N_1$ is the map-germ
$\phi_1\circ \psi_1:(\mathbb R^2,0)\rightarrow (\mathbb R^{2},0)$.
\end{defn}

\begin{thm}[\cite{Mont}]
$\mathcal K(M_1,N_1,x_1) =\mathcal K(M_2,N_2,x_2)$ if and only if the contact maps $\phi_1\circ \psi_1$ and $\phi_2\circ \psi_2$ are $\mathcal K$-equivalent.
\end{thm}

We study the contact of the surface with holomorphic curves with respect to an orthogonal complex structure $\mathcal J$ on $\mathbb R^4$. We present the results for $g_1$. The results for $g_2$ are similar. Thus we assume that $\mathcal J$ preserves the orientation of $\mathbb R^4$ and the surface $M$ is parameterized locally in Monge form at $0\in \mathbb R^4$ . It implies that the tangent space to $M$ at $0$ is spanned by $(1,0,0,0)$ and $(0,1,0,0)$ and the normal space at $0$ is spanned by $(0,0,1,0)$ and $(0,0,0,1)$. We consider holomorphic curves which are tangent to $M$ at $0$. Thus we have that $\mathcal J(1,0,0,0)=(0,1,0,0)$ and $\mathcal J(0,0,1,0)=(0,0,0,1)$. Namely, $\mathcal J|_{T_0\mathbb R^4}=\mathcal J_0|_{T_0\mathbb R^4}$. Since $\mathcal J$ is an orthogonal complex structure on $\mathbb R^4$, $\mathcal J$ is constant. Thus we have $\mathcal J=\mathcal J_0$. 
 Hence we may apply the standard identification of  $\mathbb R^4$ with
$\mathbb C^2$ given by 
$$
\mathbb R^4\ni (x_1,x_2,x_3,x_4) \mapsto (z,w)=(x_1+i x_2,x_3+ix_4)\in \mathbb C^2.
$$

 We may assume that
 a holomorphic curve in $\mathbb C^2$ is the graph of a holomorphic
function $w=f(z)$, since  $M$ is parameterized locally in Monge form at the origin. 
Thus, we analyze the contact map $k$ of $M$ in this parameterization with a holomorphic curve at the origin, and we determine 
the coefficients of degree less than $s+1$ of such holomorphic curve that has the highest order of contact in $J^s(\mathbb R^2,\mathbb R^2)$ with $M$ for $s=2, 3, 4$. Thus, 

\begin{equation}\label{cfh}
	k(x_1,x_2)=\left(\begin{array}{c}
		a(x_1,x_2)-{\text{Re}}f(x_1+ix_2)\\
		b(x_1,x_2)-{\text{Im}}f(x_1+ix_2)
	\end{array} \right).
\end{equation}

\noindent The $4$-jet of this map at $0$ has the following form
\begin{equation}\label{cfh}
j^4_0k(x_1,x_2)=\left(\begin{array}{c}
	\sum_{i+j=2}^4 \frac{a_{ij}}{i!j!}x_1^i x_2^j-{\text{Re}}j_0^4f(x_1+ix_2)\\
	 \sum_{i+j=2}^4 \frac{b_{ij}}{i!j!}x_1^i x_2^j-{\text{Im}}j_0^4f(x_1+ix_2)
	 \end{array} \right),
\end{equation}
where 
\begin{eqnarray*}
{\text{Re}}j_0^4f(x_1+ix_2)&=&
\alpha_1 x_1^2 - 2 \alpha_2 x_1 x_2 - \alpha_1 x_2^2 + \beta_1 x_1^3- 3 \beta_2 x_1^2 x_2 - 3 \beta_1 x_1 x_2^2+ \beta_2 x_2^3\\
&&+ \gamma_1 x_1^4  - 
 4 \gamma_2 x_1^3 x_2  - 6 \gamma_1 x_1^2 x_2^2  + 
 4 \gamma_2 x_1 x_2^3 + \gamma_1 x_2^4,
\end{eqnarray*}
\begin{eqnarray*}
{\text{Im}}j_0^4f(x_1+ix_2)&=&
\alpha_2 x_1^2 + 2 \alpha_1 x_1 x_2- \alpha_2 x_2^2 + \beta_2 x_1^3  + 3 \beta_1 x_1^2 x_2 - 3 \beta_2 x_1 x_2^2- \beta_1 x_2^3\\&&+ \gamma_2 x_1^4  + 
 4 \gamma_1 x_1^3 x_2 - 6 \gamma_2 x_1^2 x_2^2  - 
 4 \gamma_1 x_1 x_2^3 + \gamma_2 x_2^4.
\end{eqnarray*}

In these expressions, we use the following notation 
\begin{eqnarray*}
	a_{ij}=\frac{\partial^{i+j} a}{\partial u^i \partial v^j}(0,0)\ \  {\rm and}\ \ b_{ij}=\frac{\partial^{i+j} b}{\partial u^i \partial v^j}(0,0).
\end{eqnarray*}

Let $\mathcal E_{2,2}$ denote a $\mathcal E_2$-module of smooth map-germs $(\mathbb R^2,0)\rightarrow \mathbb R^2$.
Let $Jac(k)$ be the Jacobian module of $k$, i.e., $Jac(k)=\mathcal E_2 \left\{\frac{\partial k}{\partial x_1}, \frac{\partial k}{\partial x_2}\right\}$ and $I_k$ be the ideal $(k_1,k_2)$ in $\mathcal E_2$ generated by componets of $k$. Then $Jac(k)$ and $I_k \cdot \mathcal E_{2,2}$ are submodules of $\mathcal E_{2,2}$. The $\mathcal K$-tangent space to $k$ is the $\mathcal E_2$-module $T_k=Jac(k)+I_k \cdot \mathcal E_{2,2}$.  Since $k(0)=\frac{\partial k}{\partial x_1}(0)=\frac{\partial k}{\partial x_2}(0)=0$, the space $j^s_0(T_k)$ is a subspace of the vector space $J^s_{0,0}(2,2)$ of $s$-jets of map-germs $(\mathbb R^2,0)\rightarrow (\mathbb R^2,0)$. It is easy to see that the dimension of $J^s_{0,0}(2,2)$ equals to $s^2+3s$. We define the $s$-jet $\mathcal K$-codimension of $k$ as a codimension of $j^s_0(T_k)$ in $J_{0,0}^s(2, 2)$. A holomorphic curve has the highest possible contact with a surface $M$ in $s$-jets  at $0$ if the $s$-jet $\mathcal K$-codimension of its contact map-germ $k$ is maximal among all holomorphic curves passing throught the point $0$.

We determine the coefficients $\alpha_1$, $\alpha_2$ 
to get that the quadratic form
\begin{equation}\label{qform}
	\left(\begin{array}{ccc} \frac{1}{2}(a_{20}-2\alpha_1) &a_{11}+2\alpha_2 & \frac{1}{2}(a_{02}+2\alpha_1) \\ 
	\frac{1}{2}(b_{20}-2\alpha_2) &b_{11}-2\alpha_1 & \frac{1}{2}(b_{02}+2\alpha_2) 
	\end{array}          \right)
\end{equation}
has rank 1. We assume that the mean curvature vector $H$ of $M$ at the origin does not vanish. In this parameterization $|H|^2$ has the following form 
\begin{equation}\label{H}
|H|^2=(a_{20}+a_{02})^2+(b_{20}+b_{02})^2.
\end{equation} 

We determine coefficients $\alpha_1, \alpha_2$ such that the rank of the matrix (\ref{qform}) is not greater than $1$. By direct computation we obtain that
\begin{eqnarray*}
\alpha_1&=&\dfrac{b_{11}(a_{20}+ a_{02})^2+(a_{20}(b_{02}-a_{11})- a_{02}(a_{11}+b_{20}))(b_{20}+b_{02})}{2\left((a_{20}+a_{02})^2+(b_{20}+b_{02})^2\right)}\\
\alpha_2&=&\dfrac{-a_{11}(b_{20}+b_{02})^2+(b_{20}(a_{02}+b_{11})+b_{02}(b_{11}-a_{20}))(a_{20}+a_{02})}{2\left((a_{20}+a_{02})^2+(b_{20}+b_{02})^2\right)}
\end{eqnarray*}
and the second jet of the contact map at the origin has the following form 
\begin{equation}\label{quadraticpart}
j^2_0k(x_1,x_2)=\dfrac{Ax^2_1+2Bx_1x_2+Cx^2_2}{2\left((a_{20}+a_{02})^2+(b_{20}+b_{02})^2\right)}
\left(\begin{array}{c}
	a_{20}+a_{02}\\
	b_{20}+b_{02}
	\end{array}\right),
\end{equation}
where
\begin{eqnarray*}
A&=&(a_{20}+a_{02})(a_{20}-b_{11})+(a_{11}+b_{20})(b_{20}+b_{02}),\\
B&=&a_{20}(a_{11}-b_{02})+b_{11}(b_{20}+b_{02})+a_{02}(a_{11}+b_{20}),\\
C&=&(a_{20}+a_{02})(a_{02}+b_{11})-(a_{11}-b_{02})(b_{20}+b_{02}).	
\end{eqnarray*}
\begin{prop}\label{disc}
The discriminant $\delta$ of the quadratic form in $(\ref{quadraticpart})$ of the contact map $k$ at the origin  has the following form 
\begin{eqnarray}
\delta= -(K+K^N)|H|^2.
\end{eqnarray} 
\end{prop}
\begin{proof}
We describe the curvatures and the discriminant when the surface is parameterized locally in Monge form at the origin as follows:
\begin{eqnarray}
K=-a^2_{11}+a_{20}a_{02}-b^2_{11}+b_{20}b_{02},\ \
K^{N}=b_{11}(a_{20}-a_{02})-a_{11}(b_{20}-b_{02}),\nonumber
\end{eqnarray}
\begin{equation*}
\delta=-\left((a_{20}-b_{11})(a_{02}+b_{11})-(a_{11}-b_{02})(a_{11}+b_{20})\right)\left((a_{20}+a_{02})^2+(b_{20}+b_{02})^2\right).
\end{equation*}
Thus, using the expression of $|H|^2$ given by (\ref{H}) a direct substitution completes the proof.
\end{proof}

\begin{defn}
We say that a point $p \in M$ is $g_1$-elliptic if $(K+K^{N})(p)>0$.
We say that a point $p \in M$ is $g_1$-hyperbolic if $(K+K^{N})(p)<0$.
\end{defn}

We provide a description of the first jets of the contact map $k$ using the normal forms under the action of the contact group of local diffeomorphisms.

\begin{thm}[Classification of regular and generic singular points of $g_1$]\label{NormalForms}
Let $k$ be a contact map of $M$ and a holomorphic curve  with the highest  contact in $s$-jets at $p$ for $s=2$ in (1)-(3) and (6), $s=3$ in (4) and $s=4$ in (5).
\begin{enumerate}
\item \label{elip-form}A point $p$ is $g_1$-elliptic if and only if the contact map-germ $k$ at $p$  is $\mathcal K$-equivalent to a map-germ at $0$  with the $2$-jet of the form
\begin{equation*}
		(x^2_1+ x^2_2,\, 0).
\end{equation*}
\item A point $p$ is  $g_1$-hyperbolic if  and only if the contact map-germ $k$ at $p$  is $\mathcal K$-equivalent to a map-germ at $0$  with the $2$-jet of the form
\begin{equation*}\label{hyp-form}
		(x_1 x_2,\, 0).
\end{equation*}
\item \label{cond_sing} A point $p$ is a corank-one singular point of $g_1$ if  and only if the contact map-germ $k$ at $p$  is $\mathcal K$-equivalent to a map-germ at $0$  with the $2$-jet of the form
\begin{equation*}\label{sing-form}
		(x_1^2,\, 0).
\end{equation*}

\item \label{fold_cond} A point $p$ is a fold point of $g_1$ if and only if the contact map-germ $k$ at $p$  is $\mathcal K$-equivalent to a map-germ at $0$  with the $3$-jet of the form
\begin{equation*}\label{fold-form}
		(x_1^2+x_2^3,0).
\end{equation*}

\item \label{cusp_cond} A point $p$ is a cusp point of $g_1$ if and only if the contact map-germ $k$ at $p$  is $\mathcal K$-equivalent to a map-germ at $0$  with the $4$-jet of the form
\begin{equation*}\label{cusp-form}
		(x_1^2\pm x_2^4,\ 0).
\end{equation*}

\item \label{cond_corank2} A point $p$ is a corank-two singular point of $g_1$ if  and only if the contact map-germ $k$ at $p$  is $\mathcal K$-equivalent to a map-germ at $0$  with the $2$-jet of the form
\begin{equation*}\label{corank2-form}
		(0,\, 0).
\end{equation*}
\end{enumerate}	
\end{thm}

\begin{proof} We use the following notation $x=(x_1,x_2)$ and $k(x)=(k_1(x),k_2(x))$. Let us assume that $p=0$. First, we consider the case $|H(0)|^2\ne 0$. Then, without loss of generality, we suppose that 
\begin{equation}\label{a20+a02ne0}
a_{20}+a_{02}\ne0.
\end{equation} 
We apply the transformation 
$$k(x)\mapsto
(k_1(x), k_2(x)-\frac{b_{20}+b_{02}}{2((a_{20}+a_{02})^2+(b_{20}+b_{02})^2}k_1(x))$$ to get 
$j^2_0k(x)=(j^2_0k_1(x),0)$.
Then, Proposition \ref{disc} and classical results on quadratic forms imply (1) and (2).  In both cases the $2$-jet $\mathcal K$-codimension equals to $4$.

Now, we assume that $0$ is a corank-one singular point, i.e.,  the corank of the matrix 
\begin{equation}\label{corank-one}
{dg_1}_0=\left(\begin{array}{cc}
a_{20}-b_{11}&a_{11}-b_{02}\\
a_{11}+b_{20}&a_{02}+b_{11}
\end{array}\right)
\end{equation}
is $1$. 
It implies that
\begin{equation}\label{K+KN}
\det {dg_1}_0=(K+K^N)(0)=(a_{20}-b_{11})(a_{02}+b_{11})-(a_{11}-b_{02})(a_{11}+b_{20})=0.
\end{equation}

By (\ref{corank-one}), we may assume without loss of generality that 
\begin{equation}\label{a20neb11}
a_{20}-b_{11}\ne 0.
\end{equation}
Then, we solve the equation (\ref{K+KN}) with respect to $a_{02}$. Replacing $a_{02}$ with this solution  in $k$, we have
$$
j^2_0k(x)=\left(
\left(1 + \left(\frac{a_{11} - b_{02}}{a_{20} - b_{11}}\right)^2\right)\left((a_{20} - b_{11})x_1 + (a_{11} - b_{02})x_2\right)^2, 0 \right).
$$
Thus, we change coordinates in the following way 
$$
(x_1,x_2)\mapsto \left(\frac{x_1-(a_{11}-b_{20})x_2}{a_{20}-b_{11}},x_2\right)
$$
to get that the $2$-jet of $k$ has the following form 
$$
j^2_0k(x)=(\frac{(a_{20} - b_{11})^2 + (a_{11} + b_{20})^2}{(a_{20} - b_{11})^2}x_1^2,0). 
$$
Hence, we apply the transformation 
$$
k(x)\mapsto\left(\frac{(a_{20} - b_{11})^2}{(a_{20} - b_{11})^2 + (a_{11} + b_{20})^2}k_1(x),k_2(x)\right)
$$
to obtain that $j^2_0k(x)=(x_1^2,0)$. Its $2$-jet $\mathcal K$-codimension equals to $6$. We point out that this codimension is maximal by the construction of the normal form. Thus, we finish the proof of (\ref{cond_sing}).

The $3$-jet of $k$ at $0$ has the following form
$$
j^3_0k(x)=\left(
x_1^2+\sum_{i=0}^3k^1_{3-i,i}x_1^{3-i} x_2^i, \sum_{i=0}^3k^2_{3-i,i}x_1^{3-i} x_2^i\right),
$$
where coefficients $k^l_{3-i,i}$ depend linearly on $\beta_1, \beta_2$ for $l=1,2$ and $i=0,1,2,3$, see $(\ref{cfh})$.
Using the term $x_1^2$ in $k_1(x_1,x_2)$, we can reduce $k^2_{30}x_1^3 + k^2_{21}x_1^2x_2$ but we are not able to do the same with $k^2_{12}x_1 x_2^2+ k^2_{03}x_2^3$. Therefore, first,
we find $\beta_1$ and $\beta_2$ such that 
\begin{equation}\label{systemd1d2}
k^2_{12}=k^2_{03}=0.
\end{equation} 
The determinant of $2\times 2$ linear system (\ref{systemd1d2}) has the following form
$$
-\frac{(3 ((a_{11} - b_{02})^2 + (a_{20} - b_{11})^2)^3 ((a_{20} - b_{11})^2 + 
(a_{11} + b_{20})^2)}{(a_{20} - b_{11})^7 (a_{20} + a_{02})^2}.
$$
By (\ref{a20+a02ne0}) and (\ref{a20neb11}) this determinant is well defined and does not vanish. 
Then, we apply Wolfram
Mathematica (\cite{Math}) to find solutions $\beta_1$ and $\beta_2$. The formulas are too long to present them. Although, we compute, assuming that 
$\beta_1$ and $\beta_2$ satisfy  (\ref{systemd1d2}), to get that
\begin{equation}\label{k03}
k^1_{03}=\frac{(a_{11} - b_{02})^2 + (a_{20} - b_{11})^2}{3 (a_{20} - b_{11}) ((a_{20} - b_{11})^2 + (a_{11} + b_{20})^2))}{d(K+K^N)}_0(X(0)),
\end{equation}
where $X$ is the germ of a smooth vector field which spans  $\ker dg_1$ on $\{K+K^N=0\}$. 
By (\ref{a20neb11}), we conclude that $0$ is a fold point of $g_1$ iff $k^1_{03}$ does not vanish.

Now, we apply the transformation 
$$k(x)\mapsto(k_1(x),k_2(x)-(k^2_{30}x_1+k^2_{21}x_2)k_1(x))$$ to obtain  
$
j^3_0k(x)=\left(
x_1^2+k^1_{30}x_1^3+k^1_{21}x_1^2 x_2+k^1_{12} x_1 x_2^2+k^1_{03} x_2^3, 0
 \right)
$.
So, we assume that $k^1_{03}\ne0$ and we change coordinates in the following way 
$$
(x_1,x_2)\mapsto \left(x_1,\frac{1}{\sqrt[3]{k^1_{03}}} x_2-\frac{k^1_{12}}{3k^1_{03}}x_1\right)
$$
to express $j^3_0k$ in the following form 
$
j^3_0k(x)=(
x_1^2+k^1_{30}x_1^3+k^1_{21} x_1^2 x_2+x_2^3, 0).
$
 Then, we apply the transformation 
\begin{equation}\label{divby}
k(x)\mapsto\left(\frac{k_1(x)}{1+k^1_{30}x_1+k^1_{21}x_2},k_2(x)\right)
\end{equation} 
to get that
$
j^3_0k(x_1,x_2)=\left(x_1^2+x_2^3,0\right)
$.
Thus, the $3$-jet $\mathcal K$-codimension of $ \left(x_1^2+x_2^3,0\right)$ equals to $8$. We point out that this codimension is maximal by the construction of the normal form. Thus, we finish the proof of (\ref{fold_cond}). 

Now, we assume that $k^1_{03}=0$. By (\ref{k03}), this assumption is equivalent to 
$$d(K+K^N)_0(X(0))=0,$$ which has the following form 
\begin{align}
(a_{20} - b_{11})^3 (a_{03} + b_{12}) - (a_{20} - b_{11})^2 (a_{12} - b_{03}) (a_{11} +  b_{20}) +& \label{k03eq0} \\
 2(a_{20} - b_{11}) (a_{11} - b_{02}) (a_{21} - b_{12}) (a_{11} + 
      b_{20}) -2 (a_{20} - b_{11})^2 (a_{11} - b_{02}) (a_{12} + b_{21}) +& \nonumber \\
 (a_{20} - b_{11}) (a_{11} - b_{02})^2 (a_{21} + b_{30}) - (a_{11} - b_{02})^2 (a_{11} + b_{20}) (a_{30} - b_{21}) 
=0.& \nonumber
\end{align}
Then, we solve the equation (\ref{k03eq0}) with respect to $a_{03}$. Replacing $a_{03}$ with this solution  in $k$, we have that
$
j^3_0k(x)=\left(
x_1^2+k^1_{30}x_1^3+k^1_{21} x_1^2 x_2+k^1_{12}x_1x_2^2, 0 \right)
$.
So, we change coordinates in the following way 
$(x_1,x_2)\mapsto \left(x_1-(1/2)k^1_{12}x_2^2,x_2\right)$
to get that  
$j^3_0k(x)=\left(x_1^2+k^1_{30}x_1^3+k^1_{21} x_1^2 x_2, 0 \right)$ in the new coordinates.
Moreover, by applying transformation (\ref{divby}), we reduce $k$ to such form that 

$$
j^4_0k(x)=\left(
x_1^2+\sum_{i=0}^4k^1_{4-i,i}x_1^{4-i} x_2^i, \sum_{i=0}^4k^2_{4-i,i}x_1^{4-i} x_2^i\right),
$$
where coefficients $k^l_{4-i,i}$ depend linearly on $\beta_1, \beta_2$ for $l=1,2$ and $i=0,\cdots,4$, see (\ref{cfh}).
Now, we proceed in the same way as before.
By using the term $x_1^2$ in $k_1(x_1,x_2)$, we can reduce $k^2_{40}x_1^4 + k^2_{31}x_1^3x_2+k^2_{22}x_1^2x_2^2$ but we are not able to 
do the same with $k^2_{13}x_1 x_2^3+ k^2_{04}x_2^4$. Therefore, first,
we determine $\gamma_1$ and $\gamma_2$ such that 
\begin{equation}\label{systeme1e2}
k^2_{13}=k^2_{04}=0.
\end{equation} 
The determinant of $2\times 2$ linear system (\ref{systeme1e2}) has the following form

$$
\frac{4 ((a_{11} - b_{02})^2 + (a_{20} - b_{11})^2)^4 ((a_{20} - b_{11})^2 + (a_{11} +  b_{20})^2)}{(a_{20} - b_{11})^9 (a_{20} + a_{02})^2}.
$$
By (\ref{a20+a02ne0}) and (\ref{a20neb11}), this determinant is well defined and does not vanish. Then,  
we apply Wolfram
Mathematica (\cite{Math}) to find solutions $\gamma_1$ and $\gamma_2$. The formulas are too long to present them. 
Although, we compute, assuming that 
$\gamma_1$ and $\gamma_2$ satisfy  $(\ref{systeme1e2})$, to get that
\begin{equation}\label{k04}
k^1_{04}=\frac{(a_{11} - b_{02})^2 + (a_{20} - b_{11})^2}{12 (a_{20} - b_{11})^2 ((a_{20} - b_{11})^2 + (a_{11} + b_{20})^2)}d(d(K+K^N)(X))_0(X(0)).
\end{equation}
By (\ref{a20neb11}), we obtain that $0$ is a cusp point of $g_1$  iff $k^1_{04}$ does not vanish for $\gamma_1$, $\gamma_2$ satisfying 
$(\ref{systeme1e2})$. We assume that $k^1_{04}\ne0$.
Let us notice that the sign of $k^1_{04}$ depends only on the sign of $d(d(K+K^N)(X))_0(X(0))$ since $\frac{(a_{11} - b_{02})^2 + (a_{20} - b_{11})^2}{12 (a_{20} - b_{11})^2 ((a_{20} - b_{11})^2 + (a_{11} + b_{20})^2)}$ is positive.
We apply the transformation 
$$
k(x)\mapsto
\left(k_1(x),k_2(x)-(k^2_{40}x_1^2+k^2_{31}x_1x_2+k^2_{22}x_2^2)k_1(x)\right),
$$
to obtain that
$j^4_0k(x)=\left(
x_1^2+k^1_{40}x_1^4+k^1_{31}x_1^3 x_2+k^1_{22} x_1^2 x_2^2+k^1_{13} x_1x_2^3+k^1_{04}x_2^4, 0
 \right)$.
Then, we change coordinates in the following way 
$$
(x_1,x_2)\mapsto \left(x_1,\frac{1}{\sqrt[4]{|k^1_{04}|}}x_2-\frac{k^1_{13}}{4k^1_{04}}x_1\right)
$$
to express the $4$-jet of $k$ at $0$ in the following form 
$$
j^4_0k(x)=\left(
x_1^2+k^1_{40}x_1^4+k^1_{31} x_1^3 x_2+k^1_{22} x_1^2 x_2^2+\text{sgn}(k^1_{04}) x_2^4, 0 \right),
$$
in new coordinates.
Thus, we apply the transformation 
\begin{equation}\label{divby2}
k(x)\mapsto\left(\frac{k_1(x)}{1+k^1_{40}x_1^2+k^1_{31}x_1x_2+k^1_{22}x_2^2},k_2(x)\right)
\end{equation} 
to get that 
$j^4_0k(x_1,x_2)=\left(x_1^2\pm x_2^4,0\right)$.
Its $4$-jet $\mathcal K$-codimension equals to $10$. Thus, we finish the proof of (\ref{cusp_cond}). 

Now, we assume that $0$ is a corank-two singular point. Then, the following equations hold 
\begin{equation}\label{corank-2}
a_{20}=b_{11}=-a_{02}, \ \ b_{02}=a_{11}=-b_{20}.
\end{equation}
These conditions imply that $|H(0)|=0$, and so, the quadratic form of $g_1$ has following form 
\begin{equation}\label{qform-corank2}
	\left(\begin{array}{ccc} \frac{1}{2}(a_{20}-2\alpha_1) &-(b_{20}-2\alpha_2) & -\frac{1}{2}(a_{20}-2\alpha_1) \\ 
	\frac{1}{2}(b_{20}-2\alpha_2) &a_{20}-2\alpha_1 &-\frac{1}{2}(b_{20}-2\alpha_2). 
	\end{array}          \right).
\end{equation}
Hence, we take $\alpha_1=\frac{1}{2}a_{20}$ and $\alpha_2=\frac{1}{2}b_{20}$ to obtain (\ref{corank2-form}). Its $2$-jet $\mathcal K$-codimension equals to $10$.

Finally, we consider the case where the mean curvature vanishes at $p$. 
This condition implies that the 2-jet of the map $F$ satisfies:
\begin{eqnarray}\label{mean_vanishes}
a_{20}+a_{02}= 0, \ \ 
b_{20}+b_{02}= 0.
\end{eqnarray}
Thus we get that the quadratic part of $k$ has the following form (see (\ref{qform}))
\begin{equation}\label{qform-mean-vanishes}
	\left(\begin{array}{ccc} \frac{1}{2}(a_{20}-2\alpha_1) &a_{11}+2\alpha_2 & -\frac{1}{2}(a_{20}-2\alpha_1) \\ 
	\frac{1}{2}(b_{20}-2\alpha_2) &b_{11}-2\alpha_1 & -\frac{1}{2}(b_{20}-2\alpha_2) 
	\end{array}          \right).
\end{equation}
The discriminants of the the first and the second components are non-negative
$$(a_{11}+2\alpha_2)^2+(a_{20}-2\alpha_1)^2\ge 0,\ \ \ (b_{11}-2\alpha_1)^2+(b_{20}-2\alpha_2)^2\ge 0.$$
We assume that at least one of them is positive, since they both vanish if and only if $a_{11}+b_{20}=a_{20}-b_{11}=0$, which together with 
(\ref{mean_vanishes}) means that $0$ is a corank-two singular point. This case has already been considered. Let us choose $\alpha_1=\frac{1}{2}b_{11}$ and $\alpha_2=\frac{1}{2}b_{20}$. Then (\ref{qform-mean-vanishes}) has the following form 
\begin{equation}\label{qform-mean-vanishes-2}
	\left(\begin{array}{ccc} \frac{1}{2}(a_{20}-b_{11}) &a_{11}+b_{20} & -\frac{1}{2}(a_{20}-b_{11}) \\ 
	0 &0 &0 
	\end{array}          \right).
\end{equation}
Since $(a_{11}+b_{20})^2+(a_{20}-b_{11})^2>0$, we get that the point is $g_1$-hyperbolic and $j^2_0k$ is $\mathcal K$-equivalent to (\ref{hyp-form}). 
\end{proof}

\section{The Gauss-Bonnet type formulas}

The cusp point $p$ of $g_i$ is called positive (negative, respectively) if ${dg_i}_q: T_qM\rightarrow T_{g_i(q)}S_i$ preserves (reverses, respectively) the orientation for every $q$ in a small neighborhood $U$ of $p$ such that  $g_i|_U$ is injective at $q$.
We parameterize the singular set $\Sigma_{g_i}=\{p\in M| (K-(-1)^i K^N)(p)=0\}$ by $\xi_i$ such that $g_i(M)$ lies on the left hand side $g_i\circ \xi_i$ and let $\tau$ be the arc length parameter of  $g_i\circ \xi_i$ on $S_i$.  Then the geodesic curvature $\kappa_g$ of $g_i\circ \xi_i$ on $S_i$ is well defined in fold points since the kernel of $dg_i$ is transverse to $\Sigma_{g_i}$ at these points.
Let $dA$ be the area form on $M$ induced by the immersion. For $N\subset M$, $\chi(N)$ denotes the Euler characteristic of $N$. 
\begin{thm}[Gauss-Bonnet formulas] \label{G-B} Let $M$ be an 
oriented, closed surface.
If $M\rightarrow \mathbb R^{4}$ is an  immersion that has property $(G_2)$, then
  \begin{eqnarray}
2\pi\chi(M)&=&\int_{M}|K\pm K^N|dA+2\int_{\{K\pm K^N=0\}}\kappa_g(\tau)d\tau, \label{GB1}  \\
\frac{1}{\pi}\int_M KdA &=& \sum_{i=1}^2\chi(M_i^+)-\chi(M_i^-)+{\mathcal S}^{+}_{g_i}-{\mathcal S}^{-}_{g_i}, \label{KdA}\\
\frac{1}{\pi}\int_M K^NdA &=& \sum_{i=1}^2(-1)^{i-1}\left(\chi(M_i^+)-\chi(M_i^-)+{\mathcal S}^{+}_{g_i}-{\mathcal S}^{-}_{g_i}\right), \label{NdA}
\end{eqnarray}
where 
$$
M_i^+=\{p\in M| K(p)>(-1)^{i-1} K^{N}(p)\}, \ \ \ M_i^-=\{p\in M| K(p)<(-1)^{i-1} K^{N}(p)\},
$$
and ${\mathcal S}^+_{g_i}$ (${\mathcal S}^-_{g_i}$, respectively) is the number of positive (negative, respectively) 
cusp points of $g_i$ for $i=1,2$.
\end{thm}

\begin{proof}
By Theorem \ref{generic},  the singular sets of the components $g_i:M\rightarrow S_i$ for $i=1,2$ consist of fold points and cusp points. 
Therefore, a straightforward application of Quine's Theorem (see \cite{Q}) to $g_i$ for $i=1,2$ implies 
$$
2 \deg(g_i)=\chi(M_i^+)-\chi(M_i^-)+{\mathcal S}_{g_i}^+-{\mathcal S}_{g_i}^-,
$$ 
since $\chi(S_i)=2$.  Moreover, Theorem 4.6 in \cite{HO1} (see also Proposition 5 in \cite{We}) states
$$\int_M KdA=2\pi(\deg(g_1)+\deg(g_2)), \ \ \int_M K^NdA=2\pi(\deg(g_1)-\deg(g_2)).$$
Thus, we obtain equations (\ref{KdA})-(\ref{NdA}).

Since the singular sets of $g_1$ and $g_2$ consist of fold points and cusp points, we may apply Proposition 3.1 and Proposition 3.7 in \cite{SUY3} to conclude that 
$\int_{\Sigma_{g_i}}\kappa_g(\tau) d\tau$ is bounded. So, since the Gaussian curvature of 
$S_i$ equals to $2$ at any point, the following equation holds 
\begin{equation}\label{GB2}
2\pi \chi(M)=\int_M 2|g_i^{\ast}\sigma|+2\int_{\Sigma_{g_i}}\kappa_g(\tau) d\tau,
\end{equation}
where $\sigma$ is the area form on $S_i$. But by Proposition 4.5 in  \cite{HO1} $g_i^{\ast}\sigma=\frac{1}{2}(K-(-1)^{i} K^N)dA$,
and considering that $\Sigma_{g_i}=\{K-(-1)^{i} K^N=0\}$, we conclude that (\ref{GB2}) implies (\ref{GB1}).
\end{proof}

If the map $M\rightarrow \mathbb R^4$ is an embedding, then $\deg(g_1)=\deg(g_2)=\frac{1}{4\pi}\int_M KdA=\frac{1}{2}\chi(M)$ 
(\cite{CS}, \cite{HO1}, \cite{We}). 
Thus, as a consequence of Theorem \ref{G-B} we state.
\begin{cor}\label{G-B-1}
Let $M$ be an  
oriented, closed surface.
Let $M\rightarrow \mathbb R^{4}$ be an embedding that has property $(G_2)$. 
If 
$g: M\rightarrow S_1\times S_2$ is the Gauss map of this embedding, then for $i=1,2$ 
\begin{eqnarray}
\chi(M)= \chi(M_i^+)+\frac{1}{2}({\mathcal S}^{+}_{g_i}-{\mathcal S}^{-}_{g_i}), 
\ \ \
\chi(M_i^-)=\frac{1}{2}\left( {\mathcal S}^{+}_{g_i}-{\mathcal S}^{-}_{g_i}\right), \label{chM-} 
\end{eqnarray}
where $g_i\ i=1,2$ are the components of $g$.
\end{cor}
\begin{proof} 
Since the map $M\rightarrow \mathbb R^4$ is an embedding, we have that 
$$\frac{1}{\pi}\int_M K^NdA=0\ {\rm and}\ \frac{1}{2\pi}\int_M K dA=\chi(M)$$ 
(see \cite{We}, \cite{HO1}). Hence by (\ref{NdA}), we have that 
$$\chi(M_1^+)-\chi(M_1^-)+{\mathcal S}^{+}_{g_1}-{\mathcal S}^{-}_{g_1}=\chi(M_2^+)-\chi(M_2^-)+{\mathcal S}^{+}_{g_2}-{\mathcal S}^{-}_{g_2}.$$
 Thus (\ref{KdA}) implies that for $i=1,2$  
  \begin{equation}\label{chiM}
 \chi(M)=\chi(M_i^+)-\chi(M_i^-)+{\mathcal S}^{+}_{g_i}-{\mathcal S}^{-}_{g_i}.
 \end{equation}  
 On the other hand, we have 
 \begin{equation} \label{classic}
 \chi(M)=\chi(M_i^+)+\chi(M_i^-)+\chi(\Sigma_{g_i}).
 \end{equation}
  But $\Sigma_{g_i}$ is a disjoint union of simple, regular, closed, smooth curves. This property implies that $\chi(\Sigma_{g_i})=0$. Thus from (\ref{chiM})-(\ref{classic}), we get
  $$
  \chi(M_i^-)=\frac{1}{2}\left( {\mathcal S}^{+}_{g_i}-{\mathcal S}^{-}_{g_i}\right)\ \ {\rm and}\ \   
  \chi(M)= \chi(M_i^+)+\frac{1}{2}({\mathcal S}^{+}_{g_i}-{\mathcal S}^{-}_{g_i}),
  $$
  which finishes the proof.
\end{proof}

We finish by describing the singular set of 
a closed  
surface of genus two embedded in the standard 3-dimensional sphere of unit radius  $S^3 \subset \mathbb R^4$. 

\noindent {\bf Example 1.}


Let $T_2$ be the surface given as the 
transversal intersection of the hypersurfaces $F^{-1}(0)$ and $G_{r}^{-1}(0)$, where

$$F:\mathbb R^4 \rightarrow \mathbb R,\ \  F(u,v,x,y)=x^2- y^2 + u^3 - 3 u v^2,$$
$$\ \ \ \ \ G_{r}:\mathbb R^4 \rightarrow \mathbb R,\ \ G_{r}(u,v,x,y)=x^2+ y^2 + u^2 + v^2 - r^2.$$

\begin{figure}
\centering
	\includegraphics[width=0.4\linewidth]{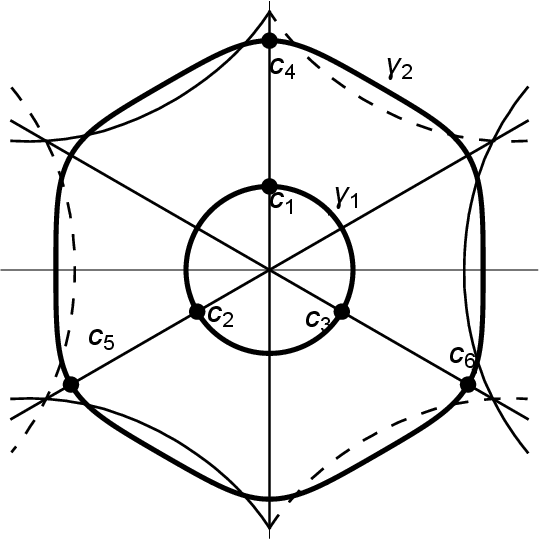}
	\caption{$h=0$ the solid curve, $k=0$ the dashed curve, the cups points $c_1,\cdots, c_6$ on the singular set $\hat \Sigma_{g_i}=\gamma_1\cup \gamma_2$.
}
	\label{fig:fig1}
\end{figure}

\noindent The projection 
\[ \pi: \mathbb{R}^{4} \rightarrow \mathbb{R}^{2},\ \ \ (u,v,x,y) \mapsto (u,v) \]

\noindent maps $T_2$ onto a hexagon $H_{uv}$ defined by the following inequalities on the $uv$-plane, see Figure \ref{fig:fig1}.

\begin{center}
$h(u,v) \geq 0,\ $  $\ k(u,v) \geq 0$, 
\end{center}
where,
\[ h(u,v) = \frac{r^{2}-u^{2}-v^{2}-u^{3}+3uv^{2}}{2},\ \ \  k(u,v) = \frac{r^{2}-u^{2}-v^{2}+u^{3}-3uv^{2}}{2}.\]

\noindent Thus, we consider four parameterizations 

\begin{eqnarray*}
{\bf x}_{\pm \pm}(u,v) = \left( u,v,\pm \sqrt{h(u,v)}, \pm \sqrt{k(u,v)} \right),
\end{eqnarray*}
with the same domain. Namely, the interior of $H_{uv} \subset \mathbb{R}^2$. 
 
We denote the image of $H_{uv}$ under ${\bf x}_{\pm,\pm}$ by $H_{\pm,\pm}$.  
This model of the double torus is obtained as the union of four identical pieces; each one 
is the closure of a hexagonal region. Namely, 

\[ T_2 = H_{++} \cup H_{+-} \cup H_{-+} \cup H_{--}, \]
see Figure \ref{fig:fig1}.
A description of the symmetries and the parameterization of this model of the double torus is presented in \cite{MVF}. Moreover, 
a direct computation implies that
the invariant 
$\Delta$, see (\ref{invariants}), is negative everywhere except at the four points parameterized by the origin in the four coordinate charts, where the Gaussian curvature is not zero. 
So,
Theorem 1.2 in \cite{L} implies that the Gauss map is a global immersion of $T_2$ into the Grassmannian. 

Now, we analyze the singularities of the components of the Gauss map.
Since the surface is spheric, the normal curvature vanishes. Thus, the singular set of both components $g_1$ and $g_2$ 
are coincident; it is the 
zero locus of the Gaussian curvature. We parameterize its projection on $H_{uv}$ by using polar coordinates as the intersection of the zero locus of the 
curve 
\begin{eqnarray*}
\hat \Sigma_{g_i}&=&9 r^{10}-248 r^8+232 r^6-464 r^4+352 r^2-32
+\left(9 r^4-86 r^2+16\right) r^6 \cos (6 \varphi ), \label{sg1}	
\end{eqnarray*}
where $i=1,2$ with the region $H_{uv}$, see Figure \ref{fig:fig1}.

By applying the symmetries of the surface, specifically, the reflections on the 
lines through the origin with angles  
\begin{equation*}
\phi=0,\ \frac{\pi}{6},\ \frac{\pi}{3},\ \frac{\pi}{2}, \frac{2\pi}{3},\frac{5\pi}{6},
\end{equation*}
is easy to get the following description.

\begin{prop}
The curve $\hat \Sigma_{g_i}$ consists of two connected components $\gamma_i, \ i=1,2$, where  $\gamma_i$ is diffeomorphic to $S^1$.
Both have index 1 for any point in the interior of the bounded region of its complement. 
\end{prop}  

We observe that the exterior component $\gamma_2$ intersects $H_{uv}$ in a slight curve near each vertex. So, each intersection 
is glued, at the boundary of the hexagon, with the other three parameterized in different coordinate charts. In this way, we get six closed curves on $T_2$.	
If we add four components parameterized by
$\gamma_1$, we conclude that the singular set of the components of the Gauss map consists of ten closed curves diffeomorphic to  
a circle. By direct computations we verify that $K<0$ in the interior of the discs bounded by these closed curves, while $K>0$ in the complement of the closure of the union of them.  Therefore, the following equations hold.
\begin{eqnarray*}
\chi(M_1^+)= \chi(M_2^+)= -12,\ {\rm and}\ \ \chi(M_1^-)=\chi(M_2^-) =10.
\end{eqnarray*}

\noindent Moreover, by the very definition of a cusp point 
we determine six cusp points as
the intersection of the singular set with the rays from the origin with angles $\phi=\frac{\pi}{2}, \frac{7\pi}{6}, \frac{9\pi}{6},$ see Figure \ref{fig:fig1}. We verify that all of them are negative, by computing of $J(g_i)$ near each cusp point.
So, there are twenty four negative cusp points on the whole surface. 
Therefore, equations  (\ref{chM-}) are satisfied in this example.

\medskip

\noindent \textbf{Acknowledgement:} 
The second and third authors were partially supported by CONACYT grant 283017 and PAPIIT-DGAPA-UNAM grant IN118217. The authors thank R\'emi Langevin for very useful disscusions.


\end{document}